\documentclass[preprint,12pt]{elsarticle}
\usepackage{geometry}
\newif\ifsubmissionpackage
\submissionpackagetrue % set to \submissionpackagefalse to obtain the clean manuscript for upload
\usepackage{lmodern}
\usepackage{anyfontsize}
\usepackage[utf8]{inputenc}
\usepackage{amsmath,amssymb,amsthm}
\usepackage{mathrsfs}
\usepackage[shortlabels]{enumitem}
\usepackage{booktabs}
\usepackage{array}
\usepackage{microtype}
\usepackage[colorlinks=true,linkcolor=green,citecolor=blue,urlcolor=green]{hyperref}
\usepackage[capitalise,nameinlink]{cleveref}
\usepackage{float}
\usepackage{orcidlink}

\numberwithin{equation}{section}
\theoremstyle{plain}
\newtheorem{theorem}{Theorem}[section]
\newtheorem{proposition}[theorem]{Proposition}
\newtheorem{lemma}[theorem]{Lemma}
\newtheorem{corollary}[theorem]{Corollary}
\theoremstyle{definition}
\newtheorem{definition}[theorem]{Definition}
\newtheorem{example}[theorem]{Example}
\theoremstyle{remark}
\newtheorem{remark}[theorem]{Remark}

\crefname{theorem}{Theorem}{Theorems}
\crefname{proposition}{Proposition}{Propositions}
\crefname{lemma}{Lemma}{Lemmas}
\crefname{corollary}{Corollary}{Corollaries}
\crefname{definition}{Definition}{Definitions}
\crefname{example}{Example}{Examples}
\crefname{remark}{Remark}{Remarks}
\crefname{section}{Section}{Sections}
\crefname{subsection}{Section}{Sections}
\crefname{appendix}{Appendix}{Appendices}
\crefname{table}{Table}{Tables}

\allowdisplaybreaks 
\newcommand{\N}{\mathbb{N}}
\newcommand{\Q}{\mathbb{Q}}
\newcommand{\R}{\mathbb{R}}
\newcommand{\cA}{\mathscr{A}}
\newcommand{\cB}{\mathscr{B}}
\newcommand{\cO}{\mathscr{O}}
\newcommand{\cC}{\mathscr{C}}

\newcommand{\cL}{\mathscr{L}}
\newcommand{\cD}{\mathscr{D}}
\newcommand{\M}{M(\cA)}
\newcommand{\1}{\mathbf{1}}
\newcommand{\cl}{\operatorname{cl}}
\newcommand{\inter}{\operatorname{int}}
\newcommand{\dist}{\operatorname{dist}}
\newcommand{\diam}{\operatorname{diam}}
\newcommand{\supp}{\operatorname{supp}}
\newcommand{\sgn}{\operatorname{sgn}}
\newcommand{\upto}{\uparrow}
\newcommand{\downto}{\downarrow}
\newcommand{\Lu}{\textup{(L)}}
\newcommand{\SK}{\textup{(SK)}}
\newcommand{\OS}{\textup{(OS)}}
\newcommand{\UA}{\textup{(UA)}}
\newcommand{\RD}{\textup{(R$\triangle$)}}
\newcommand{\Sone}{\textup{[S1]}}
\newcommand{\PS}{\textup{(S)}}
\newcommand{\PC}{\textup{(C)}}
\newcommand{\QC}{\mathrm{QC}}
\newcommand{\Lone}{L^1_C}
\newcommand{\nC}[1]{\|#1\|_{C}}
\newcommand{\nS}[1]{\|#1\|_{S}}
\newcommand{\nSh}[1]{\|#1\|_{Sh}}
\newcommand{\Ch}{(C)\!\!\int}
\newcommand{\Su}{(S)\!\!\int}
\newcommand{\Sh}{(Sh)\!\!\int}
\newcommand{\clM}{\cl_{M}}
\newcommand{\clC}{\cl_{C}}
\newcommand{\clS}{\cl_{S}}
\newcommand{\clSh}{\cl_{Sh}}
\newcommand{\Hd}{\widetilde{\mathcal{H}}^{\delta}_{\infty}}
\newcommand{\Hn}{\widetilde{\mathcal{H}}^{n}_{\infty}}
\newcommand{\Hstd}{\mathcal{H}^{\delta}_{\infty,\mathrm{std}}}

\journal{International Journal of Approximate Reasoning}

\begin{document}

\begin{frontmatter}

\title{Continuous function approximation and quasicontinuity in monotone measure spaces}

\author[a]{Thanh Ngo Vo Ke}
\ead{nvk.thanh@hutech.edu.vn}

\author[b]{Cuong Le Quoc}
\ead{cuonglq@ptit.edu.vn}

\author[c,e]{Duc Toan Truong\orcidlink{0000-0001-5562-621X}}
\ead{toan.truongduc@vlu.edu.vn}

\author[d,e]{Kiet Ho\corref{cor1}\orcidlink{0009-0003-0088-6235}}
\ead{hoanhkiet@vlu.edu.vn}

\cortext[cor1]{Corresponding authors}

\address[a]{Institute of Engineering, HUTECH University, Ho Chi Minh City, Vietnam}

\address[b]{Posts and Telecommunications Institute of Technology - Ho Chi Minh City Campus, Ho Chi Minh City 70000, Vietnam}

\address[c]{Laboratory for Chemical Computation and Modeling, Institute for Computational Science and Artificial Intelligence, Van Lang University, Ho Chi Minh City 70000, Vietnam}
\address[d]{Atomic Molecular and Optical Physics Research Group, Science and Technology Advanced Institute, Van Lang University, Ho Chi Minh City 70000, Vietnam}
\address[e]{Faculty of Applied Technology, Van Lang School of Technology, Van Lang University, Ho Chi Minh City 70000, Vietnam}

%\fntext[*]{Corresponding authors: Cuong Le Quoc, Tran Nhat Luan}
\begin{abstract}
We study the approximation of measurable functions by continuous functions, in measure and with respect to the Choquet, Sugeno and Shilkret integrals, for monotone measures on metric spaces. For every monotone measure, quasicontinuous functions are limits of continuous functions in measure and in the Sugeno integral; for unbounded Choquet integrable functions this remains true under a small-set subadditivity condition or under null-additivity. Conversely, under property [S1] and an outer smallness condition, every such limit is quasicontinuous, and examples show that neither condition can be omitted. This extends to monotone measures a classical description of the closure of continuous functions in capacitary $L^1$-spaces. We characterize the indicator functions in this closure by a symmetric approximation property of sets, and we give criteria, involving regularity and order continuity, for the density of bounded continuous functions. As applications, we show that two approximation theorems of Cao with compactly supported functions fail on spaces which are not locally compact, we prove that, for a given finite, continuous, weakly null-additive measure on a separable metric space, their conclusions for convergence in measure, for the Sugeno integral and for the Choquet integral of bounded functions hold if and only if the measure is concentrated on the set of points having a compact neighbourhood (and, under null-additivity, also for unbounded Choquet integrable functions), and we prove corrected versions on locally compact spaces; we describe the closure explicitly for possibility measures; and we give a Lebesgue-null Cantor dust which shows, for all $n\ge1$ and $0<\delta<n$, that continuous functions are not dense in the Choquet space of the dyadic Hausdorff content.
\end{abstract}


\end{frontmatter}
%\tableofcontents
%==========================================================================
\section{Introduction}\label{sec:intro}

Let $\mu$ be a finite Borel measure on a metric space $X$. Two classical facts connect measurable functions with continuous ones. Lusin's theorem says that every Borel function is continuous on a closed set whose complement has small measure, and, as a consequence, continuous functions are dense in $L^1(\mu)$. Riesz's theorem says that a sequence converging in measure (in particular, in $L^1$) has an almost everywhere convergent subsequence. Both facts have been studied intensively for non-additive set functions. Sun \cite{Sun94} showed that Riesz's theorem holds for a fuzzy measure if and only if it has property (S). Lusin's theorem was obtained by Wu and Ha \cite{WH94}, Jiang and Suzuki \cite{JS96}, Li and Yasuda \cite{LY04}, Kawabe \cite{Kaw07} and Li and Mesiar \cite{LM11}, among others, and on complete separable metric spaces by Watanabe, Yamada and Tanaka \cite{WYT13} and by Cao \cite{Cao15}. Several of these papers also contain \emph{mean approximation theorems}: a measurable function is approximated by continuous functions in measure or with respect to the Choquet or the Sugeno integral \cite{Cao15,CWL16,LM11,LY04}. We refer to the survey \cite{LMPK15} and to the monographs \cite{Den94,Pap95,WK09} for background.

A second, independent line of research comes from potential theory. For an outer, countably subadditive capacity, the closure of the continuous functions in the capacitary $L^1$-space, defined through the Choquet integral, consists exactly of the quasicontinuous functions; see \cite[Chapter~6]{AH96}, \cite{CMS11} and \cite{PS23}. Recently Harjulehto and Hurri-Syrj\"anen \cite{HH25} studied the space $\mathrm{N}L^1(\R^n,\Hd)$ for the dyadic Hausdorff content $\Hd$. They proved that a function is a limit of continuous functions in this space if and only if it is $\Hd$-quasicontinuous, that continuous functions are dense for $\delta=n$, and, in their Remark~6.3, that they are not dense for $0<\delta<n$. Related function spaces built on Hausdorff contents are studied, for instance, in \cite{BCRS25,COS24,ST22}.

The aim of this paper is to bring these two lines together. We work with an arbitrary monotone measure $\mu$ on a $\sigma$-algebra $\cA$ containing the Borel sets of a metric space, and ask the following questions. Which functions are limits of continuous functions in measure, or in the Choquet, Sugeno or Shilkret integral? When is this closure equal to the class of quasicontinuous functions? When are continuous functions dense? It turns out that a few simple conditions on $\mu$ decide these questions, and that each of them is really needed.

Our contributions may be summarized as follows.
\begin{enumerate}[(1)]
\item \emph{Riesz-type theorems in integral spaces} (\cref{sec:riesz}). Property (S) and property \Sone{} are characterized by the subsequence principle for Choquet, Sugeno and Shilkret mean convergence. These results are known in substance \cite{LMPK15,Sun94}; they serve as tools.
\item \emph{Closure theorem} (\cref{sec:closure}). For every monotone measure, quasicontinuous functions are limits of continuous functions in measure and in the Sugeno integral, and bounded ones also in the Choquet and Shilkret integrals (\cref{thm:QC-approx}(a); this is the known method of \cite{LM11,LY04}). For unbounded Choquet integrable functions we prove the same under a condition \SK{} of small-set subadditivity, which covers submeasures and relaxed subadditive measures and does not require null-additivity, and also for finite, continuous from above, null-additive measures (\cref{thm:QC-approx}(b),(c)). Conversely, property \Sone{} together with an outer smallness condition \OS{} implies that every limit in measure of continuous functions is quasicontinuous (\cref{thm:approx-QC}). Hence, under these conditions, the closure of the continuous functions is exactly the class of quasicontinuous functions (\cref{thm:closure}). Examples show that neither \OS{} nor \Sone{} can be omitted (\cref{ex:countable,ex:dense,ex:lambda}). In the case of $\Hd$ this recovers \cite[Theorem~5.11]{HH25} and shows that it is outer smallness, rather than countable subadditivity, which is decisive for the implication ``limit $\Rightarrow$ quasicontinuous''.
\item \emph{A correction} (\cref{sec:Cao}). As far as we know, it has not been observed that Theorems~4.1 and~4.2 of \cite{Cao15}, on approximation by nonnegative continuous functions with compact support, fail on every metric space which is not locally compact (\cref{ex:Cao,rem:Cao}). For a finite, continuous, weakly null-additive $\mu$ on a separable metric space we show that their conclusions for convergence in measure, for the Sugeno integral and for the Choquet integral of bounded functions hold if and only if $\mu(X\setminus U)=0$, where $U$ is the set of points having a compact neighbourhood; under null-additivity, the Choquet case for unbounded integrable $f\ge0$ follows as well (\cref{thm:Cao-char}). Consequently the two theorems, read for $f\ge0$ and with null-additivity added for unbounded Choquet integrable $f$, hold for all admissible $\mu$ if and only if $X$ is locally compact (\cref{cor:Cao-lc}); the unbounded case for weakly null-additive $\mu$ remains open (\cref{Q1}). We also prove corrected versions on locally compact separable metric spaces, and a version of \cite[Theorem~9]{LY04} on step-function approximation, for bounded functions in general and for unbounded functions under null-additivity (\cref{thm:Cao-corrected}).
\item \emph{Lusin property from regularity} (\cref{sec:lusin}). The standard scheme ``regularity $\Rightarrow$ Lusin property'' works with property \Sone{} and the pseudometric generating property (\cref{prop:regular-Lusin}).
\item \emph{Density} (\cref{sec:density}). For an arbitrary monotone measure, an indicator function $\1_A$ is a limit of continuous functions (in any of the above senses) if and only if $A$ has a symmetric approximation property \RD{} (\cref{thm:indicator}). This gives necessary and sufficient conditions for the density of bounded continuous functions in the Choquet space (\cref{thm:Choquet-density}) and in measure (\cref{thm:measure-density}); in the latter case order continuity is necessary. Under \Sone, \OS{} and \SK{}, density is equivalent to regularity and to the Lusin property (\cref{cor:density-Lusin}).
\item \emph{Hausdorff content} (\cref{sec:hausdorff}). Remark~6.3 of \cite{HH25} states that continuous functions are not dense in $\mathrm{N}L^1(\R^n,\Hd)$ for $0<\delta<n$; the examples given there cover $0<\delta\le n-1$ and $n=2$, $1<\delta<2$. We give a single example, a Lebesgue-null Cantor dust, which covers all $n\ge1$ and $0<\delta<n$ (\cref{thm:cantor}). We also give short proofs, based on our general results, of the ball example of \cite{HH25} and of the density for $\delta=n$.
\item \emph{Examples} (\cref{sec:examples}). For possibility measures the closure of the continuous functions is described explicitly, and continuous functions are dense if and only if every level set $\{\pi\ge\varepsilon\}$ of the possibility distribution is discrete (\cref{ex:possibility}). Distorted probabilities are also discussed.
\end{enumerate}
\Cref{tab:status} indicates which results are known and which appear to be new. We emphasize that the equality ``closure of continuous functions $=$ quasicontinuous functions'' is classical for outer capacities \cite{AH96,CMS11,HH25,PS23}; our contribution is its extension to monotone measures, with an identification of the conditions which make it work.

\begin{table}[H]
\centering
\small
\caption{Overview of the main results and their status.}\label{tab:status}
\begin{tabular}{@{}p{0.23\textwidth}p{0.42\textwidth}p{0.28\textwidth}@{}}
\toprule
Result & Content & Status \\
\midrule
\Cref{prop:riesz} and \Cref{prop:riesz-au} & Riesz-type theorems in Choquet/Sugeno/Shilkret mean & known in substance \cite{LMPK15,Sun94} \\
\Cref{thm:QC-approx}(a) and \;\;\; \cref{rem:known-approx} & quasicontinuous $\Rightarrow$ limit (bounded case) & known method \cite{LM11,LMPK15,LY04} \\
\Cref{thm:QC-approx}(b),(c) & unbounded Choquet case under \SK{} or null-additivity & new, as far as we know \\
\Cref{thm:approx-QC,thm:closure} & closure $=$ quasicontinuous functions under \Sone, \OS{} (and \SK) & classical for capacities; monotone-measure version new \\
\Cref{ex:countable,ex:dense,ex:lambda} & \OS{} and \Sone{} cannot be omitted & new in this context \\
\Cref{ex:Cao,rem:Cao} & \cite[Theorems~4.1--4.2]{Cao15} fail without local compactness & new (correction) \\
\Cref{lem:compact-Lusin} & increasing compact Lusin sets & known \cite{Cao15,Kaw07,WYT13}; self-contained proof \\
\Cref{thm:Cao-char,cor:Cao-lc} & Cao's conclusions (in measure, Sugeno, bounded Choquet) hold for $\mu$ iff $\mu(X\setminus U)=0$; unbounded Choquet under null-additivity; for all admissible $\mu$ iff $X$ is locally compact & new \\
\Cref{thm:Cao-corrected} & corrected theorems of Cao and of Li--Yasuda & new (correction and proof) \\
\Cref{prop:regular-Lusin} & regularity $\Rightarrow$ Lusin under \Sone{} and p.g.p. & variant of a standard scheme \cite{CWL16,JS96,Kaw07,LM11,WH94} \\
\Cref{thm:indicator,thm:Choquet-density,thm:measure-density,cor:density-Lusin} & density criteria & new \\
\Cref{prop:shilkret} & no density for the Shilkret integral & classical (weak $L^1$) \\
\Cref{thm:cantor} & Cantor dust, all $n\ge1$, $0<\delta<n$ & statement in \cite[Rem.~6.3]{HH25}; example new, as far as we know \\
\Cref{prop:ball,cor:delta-n} & ball example; density for $\delta=n$ & results of \cite{HH25}; new proofs \\
\Cref{ex:possibility} & possibility measures & new \\
\bottomrule
\end{tabular}
\end{table}

The paper is organized as follows. \Cref{sec:prelim} collects definitions and elementary facts. \Cref{sec:riesz} contains the Riesz-type theorems, \cref{sec:closure} the closure theorem, and \cref{sec:Cao} the discussion of the theorems of Cao and of Li and Yasuda. \Cref{sec:lusin} deals with regularity and the Lusin property, \cref{sec:density} with density, and \cref{sec:hausdorff} with the Hausdorff content. Further examples are given in \cref{sec:examples}, and open problems in \cref{sec:concl}. Proofs of some elementary or known auxiliary results are collected in \cref{app:proofs}.

%==========================================================================
\section{Preliminaries}\label{sec:prelim}

\subsection{Conventions}

Throughout, $(X,d)$ is a metric space, $\cO$ and $\cC$ denote the families of open and closed subsets of $X$, and $\cB$ is the Borel $\sigma$-algebra. We fix a $\sigma$-algebra $\cA$ with $\cB\subseteq\cA\subseteq 2^X$; the case $\cA=2^X$ is allowed. We write $\M$ for the set of all $\cA$-measurable \emph{finite-valued} functions $f\colon X\to\R$, $C(X)$ for the continuous functions and $C_b(X)$ for the bounded continuous functions on $X$. Continuous functions are Borel, so $f-\varphi\in\M$ whenever $f\in\M$ and $\varphi\in C(X)$. We write $A_n\upto A$ (resp.\ $A_n\downto A$) if $(A_n)$ is non-decreasing (resp.\ non-increasing) with union (resp.\ intersection) $A$.

A \emph{monotone measure} on $(X,\cA)$ is a set function $\mu\colon\cA\to[0,\infty]$ with $\mu(\emptyset)=0$ and $\mu(A)\le\mu(B)$ whenever $A\subseteq B$. (We do not require $\mu(X)>0$, as is done in \cite[Definition~2.1]{LMPK15}; this is harmless.) It is \emph{finite} if $\mu(X)<\infty$.

\noindent For $g\in\M$ with $g\ge0$ we use the Choquet, Sugeno and Shilkret integrals
\begin{align*}
 \Ch g\,d\mu&:=\int_0^\infty\mu(\{g>t\})\,dt,   \\
 \Su g\,d\mu&:=\sup_{t\ge0}\min\bigl(t,\mu(\{g\ge t\})\bigr),\\
 \Sh g\,d\mu&:=\sup_{t>0}t\,\mu(\{g\ge t\}),
\end{align*}
where the first integral is the Lebesgue integral of a non-increasing $[0,\infty]$-valued function and $0\cdot\infty:=0$. For $f\in\M$ we put $\nC f:=\Ch|f|\,d\mu$, $\nS f:=\Su|f|\,d\mu$, $\nSh f:=\Sh|f|\,d\mu$, and
\[
\Lone(\mu):=\{f\in\M:\ \nC f<\infty\}.
\]
We work with \emph{functions} and not with classes of functions equal almost everywhere: without null-additivity, $f=g$ almost everywhere does not imply $\nC f=\nC g$.
\noindent
\noindent A sequence $(f_n)\subseteq\M$ converges to $f\in\M$ \emph{in $\mu$-measure} if $\mu(\{|f_n-f|\ge\alpha\})\to0$ for every $\alpha>0$; \emph{$\mu$-almost everywhere} if there is $E\in\cA$ with $\mu(E)=0$ such that $f_n\to f$ on $X\setminus E$; and \emph{$\mu$-almost uniformly} if for every $\varepsilon>0$ there is $E\in\cA$ with $\mu(E)<\varepsilon$ such that $f_n\to f$ uniformly on $X\setminus E$ \cite[Definition~3.2]{LMPK15}. It converges \emph{in Choquet} (resp.\ \emph{Sugeno}, \emph{Shilkret}) \emph{mean} if $\nC{f_n-f}\to0$ (resp.\ $\nS{f_n-f}\to0$, $\nSh{f_n-f}\to0$).

\noindent For a family $\cD\subseteq C(X)$ we denote by $\clM(\cD)$ the set of all $f\in\M$ for which there are $\varphi_n\in\cD$ with $\varphi_n\to f$ in $\mu$-measure, and by $\clC(\cD)$, $\clS(\cD)$, $\clSh(\cD)$ the analogous sets defined by $\nC{f-\varphi_n}\to0$, $\nS{f-\varphi_n}\to0$ and $\nSh{f-\varphi_n}\to0$, respectively.

\subsection{Conditions on monotone measures}

In the following definition all sets belong to $\cA$. The first group of conditions is standard; see \cite[Definitions~2.3--2.5]{LMPK15}.

\begin{definition}\label{def:conditions}
A monotone measure $\mu$ on $(X,\cA)$ is said to be
\begin{itemize}[leftmargin=2em]
\item \emph{continuous from below} if $A_n\upto A$ implies $\mu(A_n)\to\mu(A)$; \emph{continuous from above} if $A_n\downto A$ and $\mu(A_{n_0})<\infty$ for some $n_0$ imply $\mu(A_n)\to\mu(A)$; \emph{continuous} if it is both; \emph{order continuous} if $A_n\downto\emptyset$ implies $\mu(A_n)\to0$;
\item \emph{null-additive} if $\mu(A\cup B)=\mu(A)$ whenever $\mu(B)=0$; \emph{weakly null-additive} (w.n.a.) if $\mu(A)=\mu(B)=0$ implies $\mu(A\cup B)=0$;
\item a \emph{submeasure} if $\mu(A\cup B)\le\mu(A)+\mu(B)$; \emph{countably subadditive} if $\mu(\bigcup_iA_i)\le\sum_i\mu(A_i)$; \emph{$K$-relaxed subadditive} ($K\ge1$) if $\mu(A\cup B)\le K(\mu(A)+\mu(B))$.
\end{itemize}
Moreover, $\mu$ is said to have
\begin{itemize}[leftmargin=2em]
\item the \emph{pseudometric generating property} (p.g.p.) if $\mu(A_n)\vee\mu(B_n)\to0$ implies $\mu(A_n\cup B_n)\to0$ \cite[Definition~2.5(ii)]{LMPK15};
\item \emph{property} \PS{} \cite[Definition~2.1]{Sun94} if every sequence $(A_n)$ with $\mu(A_n)\to0$ has a subsequence $(A_{n_i})$ with $\mu(\limsup_iA_{n_i})=0$;
\item \emph{property} \Sone{} \cite[Definition~5.23]{LMPK15} if every sequence $(A_n)$ with $\mu(A_n)\to0$ has a subsequence $(A_{n_i})$ with $\lim_{k\to\infty}\mu\bigl(\bigcup_{i\ge k}A_{n_i}\bigr)=0$.
\end{itemize}
Finally, we use the following three conditions.
\begin{itemize}[leftmargin=2em]
\item[\UA] For every $\varepsilon>0$ there is $\delta>0$ such that $\mu(B)<\delta$ implies $\mu(A\cup B)\le\mu(A)+\varepsilon$ for all $A$ (\emph{uniform autocontinuity from above}).
\item[\SK] There is $K\ge1$ such that for every $\varepsilon>0$ there is $\delta>0$ such that $\mu(B)<\delta$ implies $\mu(A\cup B)\le K\mu(A)+\varepsilon$ for all $A$ (\emph{small-set $K$-subadditivity}).
\item[\OS] For every $\varepsilon>0$ there is $\delta>0$ such that for every $E\in\cA$ with $\mu(E)<\delta$ there is $G\in\cO$ with $G\supseteq E$ and $\mu(G)<\varepsilon$ (\emph{outer smallness}).
\end{itemize}
\end{definition}

In \cite{LMPK15} property (S) is denoted by [S]. Relaxed subadditivity is used by Kawabe and Yamada \cite{KY23}, who require the inequality for disjoint $A,B$ only; this implies the inequality for all $A,B$, since $A\cup B=A\cup(B\setminus A)$ and $\mu(B\setminus A)\le\mu(B)$. The names \SK{} and \OS{} are ours. Condition \SK{} unifies \UA{} and relaxed subadditivity (\cref{lem:implications}(c),(d)).

\begin{remark}\label{rem:pgp}
The p.g.p.\ has the equivalent $\varepsilon$--$\delta$ form: for every $\varepsilon>0$ there is $\delta>0$ with $\mu(A\cup B)<\varepsilon$ whenever $\mu(A)\vee\mu(B)<\delta$. (If this fails for some $\varepsilon_0$, there are $A_n,B_n$ with $\mu(A_n)\vee\mu(B_n)<1/n$ and $\mu(A_n\cup B_n)\ge\varepsilon_0$, contrary to the sequential form; the converse is immediate.) By iteration, for every $\varepsilon>0$ and $m\in\N$ there is $\delta>0$ such that $\mu(A_1\cup\dots\cup A_m)<\varepsilon$ whenever $\mu(A_i)<\delta$ for all $i$. We write $\delta_{\rm pgp}(\varepsilon)$ and $\delta_{\rm OS}(\varepsilon)$ for numbers $\delta$ as in the p.g.p.\ and in \OS.
\end{remark}

\subsection{Regularity and quasicontinuity}

\begin{definition}\label{def:regular}
Let $\mu$ be a monotone measure on $(X,\cA)$ and $A\in\cA$.
\begin{enumerate}[(i)]
\item $A$ is \emph{regular} if for every $\varepsilon>0$ there are $F\in\cC$ and $G\in\cO$ with $F\subseteq A\subseteq G$ and $\mu(G\setminus F)<\varepsilon$. The measure $\mu$ is \emph{regular} if every $A\in\cA$ is regular. (This is the regularity of Wu and Ha \cite{WH94}; see also \cite[Definition~6.2]{LMPK15}.)
\item $A$ has property \RD{} (\emph{symmetric approximation}) if for every $\varepsilon>0$ there are $F\in\cC$, $G\in\cO$ with $F\subseteq G$, and $E\in\cA$ with $\mu(E)<\varepsilon$, such that $(A\triangle F)\cup(A\triangle G)\subseteq E$.
\end{enumerate}
\end{definition}

\begin{definition}\label{def:QC}
A function $f\in\M$ is \emph{$\mu$-quasicontinuous} if for every $\varepsilon>0$ there is $G\in\cO$ with $\mu(G)<\varepsilon$ such that $f|_{X\setminus G}$ is continuous; equivalently, there is $F\in\cC$ with $\mu(X\setminus F)<\varepsilon$ such that $f|_F$ is continuous. We write $\QC(\mu)$ for the set of $\mu$-quasicontinuous functions in $\M$. For $\cA=\cB$ we say that $\mu$ has the \emph{Lusin property} \Lu$[\mu]$ if every $f\in M(\cB)$ is $\mu$-quasicontinuous; this is the conclusion of the Lusin-type theorems in \cite{LM11,LY04}.
\end{definition}
\noindent For $\mu=\Hd$ this is the quasicontinuity used in \cite{HH25}.

\subsection{Elementary facts}

The following facts are elementary, and most of them are folklore; proofs are given in \cref{app:proofs} for the reader's convenience.

\begin{lemma}\label{lem:elementary}
Let $\mu$ be a monotone measure on $(X,\cA)$ and let $g,f,f_n,g_i\in\M$, $g\ge0$. Then we have
\begin{enumerate}[(a)]
\item $\int_0^\infty\mu(\{g>t\})\,dt=\int_0^\infty\mu(\{g\ge t\})\,dt$, $\sup_t\min(t,\mu(\{g>t\}))=\sup_t\min(t,\mu(\{g\ge t\}))$ and $\sup_t t\mu(\{g>t\})=\sup_t t\mu(\{g\ge t\})$.
\item \textup{(Markov)} For $\alpha>0$: $\alpha\,\mu(\{g\ge\alpha\})\le\Sh g\,d\mu\le\Ch g\,d\mu$ and $\min(\alpha,\mu(\{g\ge\alpha\}))\le\Su g\,d\mu$.
\item $\Su g\,d\mu\le\max(\alpha,\mu(\{g\ge\alpha\}))$ for every $\alpha>0$.
\item \textup{(Tail identity)} For $M\ge0$: $\Ch (g-M)^+\,d\mu=\int_M^\infty\mu(\{g>s\})\,ds$.
\item For $c\ge0$ and $A\in\cA$ we have $\Ch c\1_A\,d\mu=c\mu(A)$, $\Su c\1_A\,d\mu=\min(c,\mu(A))$ and $\Sh c\1_A\,d\mu=c\mu(A)$.
\item $\nC{cf}=|c|\,\nC f$ for $c\in\R$.
\item If $\nC f<\infty$, then $\mu(\{|f|>t\})<\infty$ for every $t>0$, $T(M):=\int_M^\infty\mu(\{|f|>s\})\,ds\to0$ as $M\to\infty$, and $\mu(\{|f|\ge t\})\le\nC f/t$.
\item If $\mu$ is $K$-relaxed subadditive, then $\mu(A_1\cup\dots\cup A_m)\le K^{m-1}\sum_i\mu(A_i)$ and $\nC{\sum_{i=1}^mg_i}\le mK^{m-1}\sum_i\nC{g_i}$; in particular $\nC{f+g}\le2K(\nC f+\nC g)$. Conversely, if $\nC{f+g}\le K'(\nC f+\nC g)$ for all $f,g\in\M$, then $\mu$ is $K'$-relaxed subadditive.
\item $\nC{f_n-f}\to0$ $\Rightarrow$ $\nSh{f_n-f}\to0$ $\Rightarrow$ $f_n\to f$ in $\mu$-measure $\Leftrightarrow$ $\nS{f_n-f}\to0$.
\end{enumerate}
\end{lemma}
\noindent Part (i) and the equivalence for the Sugeno integral are contained in \cite[Theorem~5.10(i),(ii)]{LMPK15}.

\begin{lemma}\label{lem:implications}
Let $\mu$ be a monotone measure on $(X,\cA)$. Then we have
\begin{enumerate}[(a)]
\item If $\mu$ is countably subadditive, then it has \Sone; property \Sone{} implies property \PS.
\item If $\mu$ is continuous from below and has the p.g.p., then it has \Sone. More generally, the p.g.p.\ together with property \PC{} implies \Sone, where property \PC{} of Kawabe and Yamada \cite{KY23} means: for every sequence $(E_n)\subseteq\cA$, $\sup_{l}\mu(E_k\cup\dots\cup E_{k+l})\to0$ $(k\to\infty)$ implies $\mu(\bigcup_{n\ge k}E_n)\to0$.
\item Every submeasure satisfies \UA; \UA{} implies \SK{} with $K=1$; every $K$-relaxed subadditive $\mu$ satisfies \SK; \SK{} implies the p.g.p.; and \UA{} implies null-additivity.
\item Condition \UA{} and relaxed subadditivity are independent.
\item If $\mu(E)=\inf\{\mu(G):G\in\cO,\ G\supseteq E\}$ for all $E\in\cA$, then $\mu$ has \OS{} with $\delta=\varepsilon$.
\item If every $A\in\cA$ is regular and $\mu$ has the p.g.p., then $\mu$ has \OS.
\item Every regular set has property \RD. If $\mu$ has \OS, then every set with property \RD{} is regular.
\end{enumerate}
\end{lemma}
\noindent We do not claim (b) as new: the version with property \PS{} in place of \Sone{} is \cite[Proposition~2.6(viii)]{LMPK15}, and property \PC{} was introduced in \cite{KY23} for a Cauchy criterion for convergence in measure. A short proof is included in \cref{app:proofs}. Continuity from below trivially implies property \PC.

\begin{lemma}\label{lem:indicator-QC}
Let $\mu$ be a monotone measure on $(X,\cA)$ and $A\in\cA$. Then $\1_A$ is $\mu$-quasicontinuous if and only if $A$ is regular.
\end{lemma}
\noindent 
The implication ``$\Rightarrow$'' is the elementary argument showing that the Lusin property implies regularity; we include it for completeness.
\begin{proof}
``$\Rightarrow$'': Given $\varepsilon>0$, take $F\in\cC$ with $\mu(X\setminus F)<\varepsilon$ such that $\1_A|_F$ is continuous. The sets $A\cap F=(\1_A|_F)^{-1}([1/2,\infty))$ and $F\setminus A=(\1_A|_F)^{-1}((-\infty,1/2])$ are closed in $F$, hence in $X$. Put $F':=A\cap F\subseteq A$ and $G:=X\setminus(F\setminus A)\supseteq A$; then $F'$ is closed and $G$ is open. A point $x\notin F$ lies in $G\setminus F'$, and a point $x\in F$ lies either in $F'$ or in $F\setminus A=X\setminus G$. Hence $G\setminus F'=X\setminus F$ and $\mu(G\setminus F')<\varepsilon$.\\
``$\Leftarrow$'': Let $F\subseteq A\subseteq G$ with $F\in\cC$, $G\in\cO$ and $\mu(G\setminus F)<\varepsilon$. The set $O:=G\setminus F=G\cap(X\setminus F)$ is open, and $X\setminus O=F\cup(X\setminus G)$ is the union of two disjoint closed sets, each of which is relatively open in $X\setminus O$. Since $\1_A=1$ on $F$ and $\1_A=0$ on $X\setminus G$, the function $\1_A$ is locally constant, hence continuous, on $X\setminus O$.
\end{proof}

%==========================================================================
\section{Riesz-type theorems in integral spaces}\label{sec:riesz}

The results of this section are known in substance, and we present them as tools. The Shilkret case and the explicit necessity for mean convergence are small additions.

\begin{proposition}\label{prop:riesz}
Let $\mu$ be a monotone measure on a measurable space $(X,\cA)$. The following are equivalent:
\begin{enumerate}[(i)]
\item $\mu$ has property \PS;
\item for all $f,f_n\in\M$: if $f_n\to f$ in $\mu$-measure, then some subsequence $(f_{n_k})$ converges to $f$ $\mu$-almost everywhere;
\item the same as (ii), with the hypothesis $\nC{f_n-f}\to0$;
\item the same as (ii), with the hypothesis $\nS{f_n-f}\to0$;
\item the same as (ii), with the hypothesis $\nSh{f_n-f}\to0$.
\end{enumerate}
\end{proposition}
\noindent The equivalence (i)$\Leftrightarrow$(ii) is \cite[Theorem~2.1]{Sun94} for finite-valued $\mu$ and \cite[Theorem~5.17]{LMPK15} for arbitrary monotone $\mu$ (stated there on subsets $A\in\cA$; we use $A=X$).
\begin{proof}
(i)$\Rightarrow$(ii): Choose $n_1<n_2<\cdots$ with $\mu(A_k)<1/k$, where $A_k:=\{|f_{n_k}-f|\ge1/k\}$. By \PS{} there is a subsequence $(k_i)$ with $\mu(N)=0$, $N:=\limsup_iA_{k_i}$. If $x\notin N$, then $x$ lies in only finitely many $A_{k_i}$, so $|f_{n_{k_i}}(x)-f(x)|<1/k_i$ for all large $i$.\noindent
(ii)$\Rightarrow$(iii),(iv),(v): \cref{lem:elementary}(i).\\
\noindent Each of (ii)--(v) implies (i): let $\mu(A_n)\to0$, $f_n:=\1_{A_n}$ and $f:=0$. By \cref{lem:elementary}(e), $\nC{f_n}=\nSh{f_n}=\mu(A_n)\to0$ and $\nS{f_n}=\min(1,\mu(A_n))\to0$; also $f_n\to0$ in measure. If a subsequence satisfies $\1_{A_{n_k}}\to0$ outside a set $E$ with $\mu(E)=0$, then $\limsup_kA_{n_k}\subseteq E$, because a point lying in infinitely many $A_{n_k}$ has $\1_{A_{n_k}}(x)=1$ infinitely often. Hence $\mu(\limsup_kA_{n_k})=0$.
\end{proof}

\begin{proposition}\label{prop:riesz-au}
Let $\mu$ be a monotone measure on $(X,\cA)$. The following are equivalent: \textup{(i)} $\mu$ has property \Sone; \textup{(ii)} if $f_n\to f$ in $\mu$-measure, then some subsequence converges to $f$ $\mu$-almost uniformly; \textup{(iii)--(v)} the same under convergence in Choquet, Sugeno or Shilkret mean, respectively.
\end{proposition}
\noindent The equivalence (i)$\Leftrightarrow$(ii) is \cite[Theorem~5.24(i)]{LMPK15}; see also the references given there.

\begin{proof}
(i)$\Rightarrow$(ii): With $A_k$ as in the previous proof, \Sone{} gives $(k_i)$ with $\mu(B_j)\to0$, where $B_j:=\bigcup_{i\ge j}A_{k_i}$. On $X\setminus B_j$ we have $|f_{n_{k_i}}-f|<1/k_i$ for all $i\ge j$, which is uniform convergence. (ii)$\Rightarrow$(iii)--(v): \cref{lem:elementary}(i). (iii)--(v)$\Rightarrow$(i): take $f_n=\1_{A_n}$ and $f=0$ as before. Given $\varepsilon>0$, almost uniform convergence of $\1_{A_{n_k}}$ to $0$ gives $E$ with $\mu(E)<\varepsilon$ and $j$ such that $\1_{A_{n_k}}<1$ on $X\setminus E$ for $k\ge j$, that is, $\bigcup_{k\ge j}A_{n_k}\subseteq E$. Since $j\mapsto\mu(\bigcup_{k\ge j}A_{n_k})$ is non-increasing, its limit is at most $\varepsilon$ for every $\varepsilon>0$.
\end{proof}

\begin{remark}\label{rem:riesz}
(a) Theorem~4.6 of \cite{HH25} (convergence in $\mathrm{N}L^1(\R^n,\Hd)$ implies $\Hd$-almost everywhere convergence of a subsequence) is the special case $\mu=\Hd$ of (i)$\Rightarrow$(iii) in \cref{prop:riesz}, since $\Hd$ is countably subadditive \cite[(H3)]{HH25} and hence has property \PS{} by \cref{lem:implications}(a). Our argument does not use the norm property employed in \cite{HH25}. (The functions in \cite{HH25} may take the values $\pm\infty$; see \cref{rem:extended})\\
(b) Kawabe and Yamada \cite{KY23} introduced property \PC{} to obtain a Cauchy criterion for convergence in measure, and studied completeness and separability of the space of measurable functions and of Choquet--Lorentz spaces. Cauchy-type versions of \cref{prop:riesz,prop:riesz-au} belong to this circle of ideas and are not treated here.\\
(c) Versions for pseudo-convergence (see \cite[Theorem~5.20]{LMPK15}) and for pairs of monotone measures \cite{LOM21} can be formulated in the same way; they contain nothing new.
\end{remark}

%==========================================================================
\section{The closure of continuous functions}\label{sec:closure}

\subsection{Quasicontinuous functions are limits of continuous functions}

\begin{theorem}\label{thm:QC-approx}
Let $\mu$ be a monotone measure on $(X,\cA)$, and let $f\in\M$ be $\mu$-quasicontinuous. Then we have
\begin{enumerate}[(a)]
\item There are $\varphi_n\in C(X)$ and closed sets $F_n$ with $\varphi_n=f$ on $F_n$ and $\mu(X\setminus F_n)<1/n$. Consequently $\mu(\{|\varphi_n-f|\ge\alpha\})<1/n$ for all $\alpha>0$ and $\nS{\varphi_n-f}\le1/n$. If $|f|\le M$, one can take $|\varphi_n|\le M$, and then $\nC{\varphi_n-f}\le2M/n$ and $\nSh{\varphi_n-f}\le2M/n$.
\item Suppose that $f\in\Lone(\mu)$ and that $\mu$ satisfies \SK{} (in particular, this holds if $\mu$ satisfies \UA, is a submeasure, or is $K$-relaxed subadditive). Then for every $\eta>0$ there is $\varphi\in C_b(X)\cap\Lone(\mu)$ with $\nC{f-\varphi}<\eta$. If $f\ge0$, one may take $\varphi\ge0$.
\item The conclusion of (b) also holds if, instead of \SK, $\mu$ is finite, continuous from above and null-additive.
\end{enumerate}
\end{theorem}

\begin{proof}
(a) Take closed sets $F_n$ with $\mu(X\setminus F_n)<1/n$ such that $f|_{F_n}$ is continuous. By Tietze's extension theorem (valid for closed subsets of metric spaces), $f|_{F_n}$ extends to some $\varphi_n\in C(X)$. If $|f|\le M$, we replace $\varphi_n$ by $(\varphi_n\wedge M)\vee(-M)$, which still equals $f$ on $F_n$. Then $\{|\varphi_n-f|>0\}\subseteq X\setminus F_n$, which gives the bound for the measure, and
\[
\nS{\varphi_n-f}=\sup_t\min\bigl(t,\mu(\{|\varphi_n-f|\ge t\})\bigr)\le\sup_{t>0}\mu(X\setminus F_n)\le1/n
\]
If $|f|,|\varphi_n|\le M$, then $|\varphi_n-f|\le2M\,\1_{X\setminus F_n}$, and \cref{lem:elementary}(e) together with the monotonicity of the integrals gives the bounds for the Choquet and Shilkret integrals.\\
(b) Let $K$ and $\delta(\cdot)$ be as in \SK{}, and let $T(M):=\int_M^\infty\mu(\{|f|>s\})\,ds$, so that $T(M)\to0$ by \cref{lem:elementary}(g). Fix $\eta>0$. We make the choices in the following order: first $M>0$ with $T(M)<\eta$; then $\varepsilon:=\eta/(2M)$; then $\delta:=\delta(\varepsilon)$; and then, by quasicontinuity, a closed set $F$ with $\mu(E)<\delta$, where $E:=X\setminus F$, such that $f|_F$ is continuous. Put $f_M:=(f\wedge M)\vee(-M)$. Then $f_M|_F$ is continuous; let $\tilde\varphi$ be a Tietze extension of $f_M|_F$ and $\varphi:=(\tilde\varphi\wedge M)\vee(-M)\in C_b(X)$. Thus $\varphi=f_M$ on $F$ and $|\varphi|\le M$.\\
On $F$ we have $|f-\varphi|=|f-f_M|=(|f|-M)^+$, so $F\cap\{|f-\varphi|>t\}\subseteq\{|f|>M+t\}$ for $t\ge0$. Off $F$ we have $|f-\varphi|\le|f|+M$. Hence
\begin{align*}
\{|f-\varphi|>t\}&\subseteq\{|f|>M+t\}\cup E &&\text{for }0\le t<2M,\\
\{|f-\varphi|>t\}&\subseteq\{|f|>M+t\}\cup\{|f|>t-M\}=\{|f|>t-M\} &&\text{for }t\ge2M.
\end{align*}
Since $\mu(E)<\delta$, \SK{} gives $\mu(\{|f|>M+t\}\cup E)\le K\mu(\{|f|>M+t\})+\varepsilon$. Therefore,
\begin{align*}
\nC{f-\varphi}&\le\int_0^{2M}\bigl[K\mu(\{|f|>M+t\})+\varepsilon\bigr]dt+\int_{2M}^\infty\mu(\{|f|>t-M\})\,dt\\
&\le K\int_M^{3M}\mu(\{|f|>s\})\,ds+2M\varepsilon+T(M)\le(K+2)\eta .
\end{align*}
Moreover $\varphi\in\Lone(\mu)$: on $F$ we have $|\varphi|=|f_M|\le|f|$, and $|\varphi|\le M$ everywhere, so $\{|\varphi|>t\}\subseteq\{|f|>t\}\cup E$ for $t<M$ and $\{|\varphi|>t\}=\emptyset$ for $t\ge M$; thus $\nC\varphi\le\int_0^M(K\mu(\{|f|>t\})+\varepsilon)\,dt\le K\nC f+M\varepsilon<\infty$. If $f\ge0$, we replace $\varphi$ by $\varphi\vee0$: then $|f-\varphi\vee0|\le|f-\varphi|$ pointwise and $\varphi\vee0=f_M$ on $F$. Since $\eta>0$ is arbitrary, (b) follows.\\
(c) Let $F_k$ be closed sets with $\mu(X\setminus F_k)<1/k$ such that $f|_{F_k}$ is continuous, and put $F'_k:=F_1\cup\dots\cup F_k$. The restriction $f|_{F'_k}$ is continuous by the pasting lemma (finitely many closed sets, on each of which $f$ is continuous), $\mu(X\setminus F'_k)<1/k$, and $E_k:=X\setminus F'_k$ decreases to a set $E_\infty$ with $\mu(E_\infty)\le\mu(E_k)<1/k$ for all $k$; so $\mu(E_\infty)=0$. Fix $\eta>0$ and $M$ with $T(M)<\eta$. For each $t\in[0,2M]$, the sets $\{|f|>M+t\}\cup E_k$ decrease in $k$ to $\{|f|>M+t\}\cup E_\infty$. Since $\mu$ is finite and continuous from above, their measures converge to $\mu(\{|f|>M+t\}\cup E_\infty)$, which equals $\mu(\{|f|>M+t\})$ by null-additivity. The functions $t\mapsto\mu(\{|f|>M+t\}\cup E_k)$ are non-increasing (hence measurable) and bounded by $\mu(X)<\infty$, so by the bounded convergence theorem
\[
\int_0^{2M}\mu(\{|f|>M+t\}\cup E_k)\,dt\longrightarrow\int_M^{3M}\mu(\{|f|>s\})\,ds\le T(M)<\eta .
\]
Choose $k$ such that the left-hand side is less than $2\eta$, and construct $\varphi$ from $F:=F'_k$ as in (b). The inclusions for the level sets obtained in (b) give $\nC{f-\varphi}<2\eta+T(M)<3\eta$. Finally $\nC\varphi\le M\mu(X)<\infty$.
\end{proof}
\noindent
Tietze's theorem is where the metric structure of $X$ is used. In (b) and (c) the only difficulty is the union $\{|f|>M+t\}\cup(X\setminus F)$. Condition \SK{} controls it uniformly. In (c), continuity from above and null-additivity control it along an increasing sequence of Lusin sets; these sets must be made increasing (by the pasting lemma), and only in (b) can the Lusin set be chosen after $M$.

\begin{remark}\label{rem:known-approx}
The mean approximation theorems \cite[Theorems~5 and~6]{LY04} (finite, continuous, weakly null-additive $\mu$ on $\cB$, Borel $f$; the Choquet part for bounded $f$), \cite[Theorems~6.13, 6.14 and 6.16]{LMPK15} and \cite[Theorem~5.4]{LM11} are obtained by combining \cref{thm:QC-approx}(a) with the corresponding Lusin theorems \cite[Theorem~4]{LY04}, \cite[Theorems~6.7 and 6.11]{LMPK15} and \cite[Corollary~5.3]{LM11}, and, for the almost everywhere approximation in \cite[Theorem~5.4]{LM11}, with \cref{prop:riesz}. This is exactly how these results are proved in the cited papers.
\end{remark}

\begin{corollary}\label{cor:unbounded}
Let $\mu$ be a finite, continuous and null-additive monotone measure on $\cB$. Then every $f\in\Lone(\mu)$ is a Choquet mean limit of bounded continuous functions, which may be taken nonnegative if $f\ge0$. The same conclusion holds if $\mu$ has the Lusin property \Lu$[\mu]$ and satisfies \SK.
\end{corollary}

\begin{proof}
A null-additive measure is weakly null-additive (see also \cite[Proposition~2.6(ii)]{LMPK15}), so $f$ is $\mu$-quasicontinuous by \cite[Theorem~4]{LY04}, which applies to finite, continuous, weakly null-additive $\mu$. Now apply \cref{thm:QC-approx}(c). The second statement follows from \cref{thm:QC-approx}(b).
\end{proof}

\subsection{Limits of continuous functions are quasicontinuous}

\begin{theorem}\label{thm:approx-QC}
Let $\mu$ be a monotone measure on $(X,\cA)$ with property \Sone{} and condition \OS. If $\varphi_n\in C(X)$ and $\varphi_n\to f\in\M$ in $\mu$-measure (in particular, if $\nC{\varphi_n-f}\to0$, $\nS{\varphi_n-f}\to0$ or $\nSh{\varphi_n-f}\to0$; see \cref{lem:elementary}(i)), then $f$ is $\mu$-quasicontinuous.
\end{theorem}

\begin{proof}
Choose $n_1<n_2<\cdots$ with $\mu(A_k)<1/k$, where $A_k:=\{|\varphi_{n_k}-f|\ge2^{-k}\}\in\cA$. By \Sone{} there are $k_1<k_2<\cdots$ with $\mu(B_j)\to0$, where $B_j:=\bigcup_{i\ge j}A_{k_i}$. Let $\varepsilon>0$ and $\delta:=\delta_{\rm OS}(\varepsilon)$; choose $j$ with $\mu(B_j)<\delta$, and $G\in\cO$ with $G\supseteq B_j$ and $\mu(G)<\varepsilon$. For $x\in X\setminus G$ and $i\ge j$ we have $|\varphi_{n_{k_i}}(x)-f(x)|<2^{-k_i}$. Thus the continuous functions $\varphi_{n_{k_i}}|_{X\setminus G}$ converge uniformly to $f|_{X\setminus G}$, which is therefore continuous.
\end{proof}

For $\mu=\Hd$ this is the argument in the proof of \cite[Theorem~5.11]{HH25}; property \Sone{} replaces the countable subadditivity (H3), and \OS{} replaces the use of the definition of the content through open covers.

\begin{remark}\label{rem:credit}
For outer, countably subadditive capacities of potential theory, the fact that the closure of the continuous functions in the capacitary $L^1$-space consists exactly of the quasicontinuous functions is classical; see \cite[Chapter~6]{AH96}, \cite{CMS11}, \cite[Proposition~3.2]{PS23} (as credited in \cite{HH25}) and \cite[Theorem~5.11]{HH25} for $\Hd$. What we add in \cref{thm:QC-approx,thm:approx-QC,thm:closure} is: the extension to arbitrary monotone measures on arbitrary $\sigma$-algebras $\cA\supseteq\cB$; the identification of the two structural hypotheses, \Sone{} for ``limit $\Rightarrow$ quasicontinuous'' and \SK{} (or null-additivity) for ``quasicontinuous $\Rightarrow$ limit'', together with \OS{} and the examples in \cref{subsec:sharp} showing that they cannot be omitted; the case of \cref{thm:QC-approx}(b) for measures which are not null-additive; and the consequences for classes of monotone measures used in fuzzy measure theory (\cref{sec:Cao,sec:examples}).
\end{remark}

\begin{remark}\label{rem:route2}
One might try to replace \Sone{} in \cref{thm:approx-QC} by property \PS{} together with condition (E) and the p.g.p., using Riesz's theorem and an Egoroff theorem with closed sets \cite[Theorem~6.10]{LMPK15}. This is correct but of no interest: under condition (E) and the p.g.p.\ every $f\in M(\cB)$ is already quasicontinuous by \cite[Theorem~6.11]{LMPK15}, which is \cite[Theorem~5.3]{LM11}.
\end{remark}

\begin{theorem}\label{thm:closure}
Let $\mu$ be a monotone measure on $(X,\cA)$. Then we have
\begin{enumerate}[(a)]
\item For every $\mu$: $\QC(\mu)\subseteq\clM(C(X))=\clS(C(X))$, and every bounded $f\in\QC(\mu)$ lies in $\clC(C_b(X))\cap\clSh(C_b(X))$.
\item If $\mu$ has \Sone{} and \OS, then $\clM(C(X))=\clS(C(X))=\QC(\mu)$ and $\clC(C(X))\cup\clSh(C(X))\subseteq\QC(\mu)$.
\item If $\mu$ has \Sone, \OS{} and \SK, then for $f\in\Lone(\mu)$ the following are equivalent: $f\in\clC(C(X))$; $f\in\clC(C_b(X)\cap\Lone(\mu))$; $f$ is $\mu$-quasicontinuous.
\end{enumerate}
\end{theorem}

\begin{proof}
(a) \Cref{thm:QC-approx}(a) and \cref{lem:elementary}(i). (b) The inclusion ``$\supseteq$'' is (a), and ``$\subseteq$'' is \cref{thm:approx-QC}. (c) \Cref{thm:QC-approx}(b) and \cref{thm:approx-QC}.
\end{proof}

\begin{corollary}[{\cite[Theorem~5.11]{HH25}} for finite-valued functions]\label{cor:HH}
Let $0<\delta\le n$, $\cA=2^{\R^n}$ and $\mu=\Hd$. Then a finite-valued $f$ with $\Ch|f|\,d\Hd<\infty$ is a limit of continuous functions with respect to $\Ch|\cdot|\,d\Hd$ if and only if it is $\Hd$-quasicontinuous.
\end{corollary}

\begin{proof}
The content $\Hd$ is countably subadditive \cite[(H3)]{HH25}; hence it has \Sone{} (\cref{lem:implications}(a)) and \SK{} with $K=1$ (\cref{lem:implications}(c)). It also has \OS{} with $\delta=\varepsilon$: if $\Hd(E)<\varepsilon$, there are dyadic cubes $Q_i$ with $E\subseteq O:=\inter(\bigcup_iQ_i)$ and $\sum_i\ell(Q_i)^\delta<\varepsilon$; the same family is admissible for the open set $O$, so $\Hd(O)<\varepsilon$. Now apply \cref{thm:closure}(c).
\end{proof}

In \cite{HH25} the implication ``quasicontinuous $\Rightarrow$ limit'' is attributed to Ponce and Spector \cite[Proposition~3.2]{PS23}, and the converse is proved as in \cref{thm:approx-QC}. \Cref{cor:HH} concerns the setting of \cite[Definition~5.8]{HH25}, where ``$\mathrm{N}L^1$-limit of continuous functions'' is defined for $f\in\mathrm{N}L^1(\R^n,\Hd)$; in \cite[Definition~4.3]{HH25} one works with $\Hd$-almost everywhere classes of $[-\infty,\infty]$-valued functions. The next remark shows that this setting reduces to ours. The examples in \cref{subsec:sharp} show that, in \cite[Theorem~5.11]{HH25}, the decisive property for the implication ``limit $\Rightarrow$ quasicontinuous'' is outer smallness, not countable subadditivity.

\begin{remark}[Extended real-valued functions]\label{rem:extended}
Let $\mu$ be null-additive with the p.g.p.\ and \OS, and let $f\colon X\to[-\infty,\infty]$ be $\cA$-measurable with $\Ch|f|\,d\mu<\infty$ (with the convention that $\{|f|>t\}$ contains $\{|f|=\infty\}$). Put $N:=\{|f|=\infty\}$, $\tilde f:=f$ on $X\setminus N$ and $\tilde f:=0$ on $N$. Then $\mu(N)\le\mu(\{|f|>t\})\le\nC f/t\to0$, so $\mu(N)=0$. For $\varphi\in C(X)$ and $t\ge0$ put $S:=\{|f-\varphi|>t\}$ and $\tilde S:=\{|\tilde f-\varphi|>t\}$ (both in $\cA$); they coincide off $N$, i.e.\ $S\setminus N=\tilde S\setminus N$. Null-additivity gives $\mu(S)\le\mu(S\cup N)=\mu((S\setminus N)\cup N)=\mu(S\setminus N)\le\mu(S)$, so $\mu(S)=\mu(S\setminus N)$, and likewise $\mu(\tilde S)=\mu(\tilde S\setminus N)=\mu(S\setminus N)$. Hence $\mu(S)=\mu(\tilde S)$ for all $t$, and $\nC{f-\varphi}=\nC{\tilde f-\varphi}$. Moreover, $f$ is quasicontinuous (with continuity into $[-\infty,\infty]$) if and only if $\tilde f$ is: given $\varepsilon>0$, cover $N$ by an open set $O'$ of small measure, using \OS, and combine it with the exceptional open set of the other function by the p.g.p. Hence the extended real-valued case reduces to ours. The content $\Hd$ satisfies all three hypotheses, since it is subadditive (\cref{lem:implications}(c)).
\end{remark}

\subsection{The hypotheses \OS{} and \Sone{} cannot be omitted}\label{subsec:sharp}

\begin{example}[\OS{} fails]\label{ex:countable}
Let $X=[0,1]$, $\cA=\cB$ (or $2^X$), and $\mu(A)=1$ if $A$ is uncountable, $\mu(A)=0$ otherwise. Then we have
\begin{itemize}[leftmargin=2em]
\item $\mu$ is monotone and countably subadditive: if $\bigcup_iA_i$ is uncountable, then some $A_i$ is uncountable. Hence $\mu$ has \Sone, \PS, the p.g.p., \UA{} and \SK{} (\cref{lem:implications}).
\item $\mu$ is continuous from below (a countable union of countable sets is countable), but not order continuous: $(0,1/n)\downto\emptyset$ and $\mu((0,1/n))=1$.
\item \OS{} fails: $E:=\Q\cap[0,1]$ has $\mu(E)=0$, but every open $G\supseteq E$ is nonempty and open, hence uncountable, so $\mu(G)=1$.
\item Let $f:=\1_{\Q\cap[0,1]}$ and $\varphi_n:\equiv0$. Then $\nC{f-0}=\mu(\Q\cap[0,1])=0$ (\cref{lem:elementary}(e)), so $f\in\clC(C_b)$, and also $f\in\clS(C_b)\cap\clSh(C_b)\cap\clM(C_b)$.
\item $f$ is not quasicontinuous: for $\varepsilon\le1$, an open $G$ with $\mu(G)<\varepsilon$ is countable, and a countable open subset of $[0,1]$ is empty; but $f$ is not continuous on $[0,1]$.
\item The set $\Q\cap[0,1]$ has property \RD{} (take $E=\Q\cap[0,1]$ and $F=G=\emptyset$) but is not regular (\cref{lem:indicator-QC}); so \OS{} is needed in \cref{lem:implications}(g). The set $[0,1/2]$ does not have property \RD: if $(A\triangle F)\cup(A\triangle G)\subseteq E$ with $E$ countable and $A=[0,1/2]$, then $F\supseteq[0,1/2]\setminus E$ is dense in $[0,1/2]$ and closed, so $F\supseteq[0,1/2]$; then the open set $G\supseteq[0,1/2]$ contains $[0,1/2+\eta)$ for some $\eta>0$, and $G\setminus A\supseteq(1/2,1/2+\eta)$ is uncountable. Hence $C_b$ is not dense in $\Lone(\mu)$ (\cref{thm:Choquet-density}).
\end{itemize}
The same phenomenon occurs for the capacity $\lambda$ of \cref{ex:lambda} and for the capacity $\nu$ defined on $[0,1]$ by $\nu(A)=1$ if $A$ is non-meagre and $\nu(A)=0$ otherwise. The latter is countably subadditive, since meagre sets form a $\sigma$-ideal, and continuous from below; $\1_\Q$ has $\nu$-integral $0$, but it is not $\nu$-quasicontinuous, because an open $G$ with $\nu(G)<1$ is meagre, hence empty by Baire's theorem.
\end{example}

\begin{example}[\Sone{} and \PS{} fail, $\mu$ outer regular]\label{ex:dense}
Let $X=[0,1]$, $\cA=\cB$, and $\mu(A)=1$ if $A$ is dense in $[0,1]$, $\mu(A)=0$ otherwise.
\begin{itemize}[leftmargin=1em]
\item $\mu$ is monotone, since supersets of dense sets are dense.
\item $\mu(E)=\inf\{\mu(G):G\in\cO,\ G\supseteq E\}$ for all $E$, hence \OS{} holds with $\delta=\varepsilon$ (\cref{lem:implications}(e)): if $A$ is dense, every open $G\supseteq A$ is dense; if $A$ is not dense, there is a nonempty relatively open interval $I\subseteq[0,1]$ with $I\cap A=\emptyset$; pick $a<b$ with $[a,b]\subseteq I$; then $G:=[0,1]\setminus[a,b]$ is open, $G\supseteq A$, and $G$ misses $(a,b)$, so $\mu(G)=0$.
\item Property \PS{} (and hence \Sone) fails: let $(q_j)$ enumerate $\Q\cap[0,1]$ and $A_n:=\{q_1,\dots,q_n\}$, which is finite and not dense, so $\mu(A_n)=0$. For every subsequence, $\limsup_iA_{n_i}=\Q\cap[0,1]$ (each $q_j$ lies in $A_{n_i}$ whenever $n_i\ge j$), which is dense.
\item $\mu$ has neither the p.g.p.\ nor weak null-additivity: $\mu(\Q\cap[0,1/2])=\mu(\Q\cap[1/2,1])=0$, but the union is dense.
\item Let $J_j:=(2^{-j},2^{1-j}]$ and $f:=\sum_{j\ge1}\frac1j\1_{\Q\cap J_j}$, a Borel function with $0\le f\le1$. For $\alpha>0$ and $K:=\lfloor1/\alpha\rfloor$ we have $\{f\ge\alpha\}\subseteq\Q\cap(2^{-K},1]$, which is not dense. Thus $\mu(\{f\ge\alpha\})=0=\mu(\{f>t\})$ for all $\alpha,t>0$, and $\varphi_n:\equiv0$ satisfies $\nC{f-0}=\nS{f-0}=\nSh{f-0}=0$.
\item $f$ is not quasicontinuous: for $\varepsilon\le1$, an open $G$ with $\mu(G)<\varepsilon$ is not dense, so $X\setminus G$ contains a nondegenerate interval $[a,b]$. For some $j$, $(a,b)\cap J_j$ contains a nondegenerate interval, on which $f$ takes the value $1/j$ on the dense set of rationals and $0$ on the dense set of irrationals; so $f|_{[a,b]}$ is not continuous.
\end{itemize}
(A simpler function such as $\1_{\Q\cap(1/2,1]}$ would not do: it is quasicontinuous, with $G=(1/3,1]$.)
\end{example}

\begin{example}[The essential supremum capacity]\label{ex:lambda}
Let $X=[0,1]$, $\cA=\cB$, $m$ the Lebesgue measure, and $\lambda(A)=1$ if $m(A)>0$, $\lambda(A)=0$ otherwise. Then $\lambda$ is countably subadditive (if $m(\bigcup_iA_i)>0$, then some $m(A_i)>0$), and \OS{} fails ($\Q\cap[0,1]$, as in \cref{ex:countable}). We have $\Ch|f|\,d\lambda=\int_0^\infty\1[m(\{|f|>t\})>0]\,dt=\operatorname{ess\,sup}|f|$, so $\Lone(\lambda)=L^\infty(m)$ as a space of functions with the essential supremum as prenorm. A function is $\lambda$-quasicontinuous if and only if it is continuous, since an open $G$ with $\lambda(G)<1$ satisfies $m(G)=0$ and is therefore empty. On the other hand,
\[
\clC(C(X))=\{f\in M(\cB):\ f=g\ m\text{-a.e.\ for some }g\in C[0,1]\}.
\]
Indeed, ``$\supseteq$'' follows with the constant sequence $g$. For ``$\subseteq$'', note that for continuous $\psi$ the set $\{|\psi|>t\}$ is open, so $\operatorname{ess\,sup}|\psi|=\sup|\psi|$; hence $\operatorname{ess\,sup}|f-\varphi_n|\to0$ implies that $(\varphi_n)$ is uniformly Cauchy, $\varphi_n\to g\in C[0,1]$ uniformly, and $\operatorname{ess\,sup}|f-g|\le\operatorname{ess\,sup}|f-\varphi_n|+\sup|\varphi_n-g|\to0$. For instance, $\1_\Q\in\clC(C(X))$ is not quasicontinuous.
\end{example}

\begin{remark}[\Sone{} is not necessary in \cref{thm:approx-QC}]\label{rem:S1-not-nec}
Let $X=[0,1]$ and $v(A)=1$ if $A=X$, $v(A)=0$ otherwise. Then $v$ fails \PS, \Sone{} and the p.g.p.\ (for $A_n=X\setminus\{q_n\}$ we have $v(A_n)=0$, while every tail union equals $X$), but it has \OS{} (if $v(E)<1$, pick $x\notin E$; then $G:=X\setminus\{x\}$ is open with $v(G)=0$), and every $f$ is $v$-quasicontinuous ($G=X\setminus\{x_0\}$ has $v(G)=0$, and $f|_{\{x_0\}}$ is continuous). So the conclusion of \cref{thm:approx-QC} holds trivially. Here $\nC f=\inf|f|$ (the set $\{|f|>t\}$ equals $X$ if and only if $t<\inf|f|$, up to the single value $t=\inf|f|$), a degenerate prenorm, and the closure of the constants in the Choquet sense is all of $\M$. The meaningful question is whether \Sone{} can be weakened to \PS{} under the p.g.p.\ and \OS; see \cref{sec:concl}.
\end{remark}

%==========================================================================
\section{Compact supports: a correction of theorems of Cao and of Li and Yasuda}\label{sec:Cao}

Let $(X,\rho)$ be a complete separable metric space and $\mu$ a finite, continuous, weakly null-additive fuzzy measure on $\cB$. Theorem~4.1 of \cite{Cao15} asserts that for every measurable $f$ there are $\psi_n\in C_0^+(X)$ (nonnegative continuous functions with compact support) with $\psi_n\to f$ in $\mu$-measure. Theorem~4.2 of \cite{Cao15} asserts that, for Choquet (resp.\ Sugeno) integrable $f$, i.e.\ $f\ge0$ with $\Ch f\,d\mu<\infty$ (resp.\ $\Su f\,d\mu<\infty$), there are $\psi_k\in C_0^+(X)$ with $\Ch|\psi_k-f|\,d\mu\to0$ (resp.\ $\Su|\psi_k-f|\,d\mu\to0$). Both theorems are stated without proof. The following example shows that they are not correct as stated.

\begin{example}\label{ex:Cao}
(i) \emph{A space which is not locally compact.} Let $X=\ell^2$ (complete and separable), $\mu=\delta_0$ the Dirac measure at $0$ on $\cB(\ell^2)$, and $f\equiv1$. Then
\begin{itemize}[leftmargin=1em]
\item $\mu$ is a $\sigma$-additive probability measure; hence it is finite, continuous from above and from below, null-additive and weakly null-additive. The function $f$ is Borel, and $\Ch f\,d\mu=\Su f\,d\mu=1$.
\item $C_0^+(\ell^2)=\{0\}$. Indeed, let $\psi\in C_0^+(\ell^2)$ with $\psi(x_0)>0$. The set $W:=\{\psi>0\}$ is open and nonempty, so it contains a closed ball $\bar B(x_0,r)$ with $r>0$. This ball is closed and contained in the compact set $\supp\psi$, so it is compact. But no closed ball of positive radius in an infinite-dimensional normed space is compact (by F.~Riesz's lemma; in $\ell^2$ the unit vectors satisfy $\|e_k-e_j\|=\sqrt2$ for $k\ne j$). This is a contradiction.
\item Hence $\psi_k\equiv0$ for all $k$, $|\psi_k-f|\equiv1$, $\delta_0(\{|\psi_k-f|\ge\sigma\})=1$ for $0<\sigma\le1$, and $\Ch|\psi_k-f|\,d\delta_0=\Su|\psi_k-f|\,d\delta_0=1$ for every $k$. Both \cite[Theorem~4.1]{Cao15} and \cite[Theorem~4.2]{Cao15} fail.
\end{itemize}
(ii) \emph{Signed functions.} Theorem~4.1 of \cite{Cao15} is stated for arbitrary measurable $f$, with nonnegative approximants. For $f\equiv-1$ and any $\mu$ with $\mu(X)>0$, every $\psi\ge0$ satisfies $|\psi-f|\ge1$, so $\mu(\{|\psi-f|\ge\sigma\})=\mu(X)$ for $\sigma\le1$, and convergence in measure fails on every space. This is probably only a misprint.
\end{example}
\noindent Compact supports require local compactness, at least on a set which carries the measure, and the unbounded Choquet case requires an additional argument, which we can supply under null-additivity. We first need compact Lusin sets.

\begin{lemma}\label{lem:compact-Lusin}
Let $X$ be a Polish space, $\mu$ a finite, continuous, weakly null-additive monotone measure on $\cB$, and $f$ a Borel function. Then there are compact sets $C_1\subseteq C_2\subseteq\cdots$ such that $f|_{C_k}$ is continuous for every $k$, $\mu(X\setminus C_k)\to0$, and $N:=\bigcap_k(X\setminus C_k)$ satisfies $\mu(N)=0$.
\end{lemma}
\noindent Tightness of finite continuous fuzzy measures on complete separable metric spaces is \cite[Lemma~3.2]{Cao15}, and it is also contained in \cite[Theorem~17 and Corollary~20]{WYT13}. Compact Lusin sets are stated in \cite[Theorem~3.3]{Cao15}, without proof. The intersection device used in step (3) of our proof is that of \cite[proof of Theorem~3]{Kaw07}, where two increasing sequences of closed sets are intersected and weak null-additivity is applied to the union of the two null residuals; for real-valued measures that proof yields increasing closed Lusin sets for weakly null-additive fuzzy measures on metric spaces. For completeness, a self-contained proof is given in \cref{app:compact-Lusin}.

We first characterize, for each individual measure, when the conclusions of \cite[Theorems~4.1--4.2]{Cao15} hold. For a metric space $X$ let
\[
U:=\{x\in X:\ x\text{ has a compact neighbourhood}\}.
\]
Thus $X$ is locally compact if and only if $U=X$. We write $C_c^+(X)$ for the nonnegative functions in $C_c(X)$; this is the class $C_0^+(X)$ of \cite{Cao15}.

\begin{theorem}\label{thm:Cao-char}
Let $X$ be a separable metric space, and $\mu$ a finite, continuous and weakly null-additive monotone measure on $\cB$. The following are equivalent:
\begin{enumerate}[(a)]
\item for every Borel function $f$ there are $\psi_n\in C_c(X)$ with $\psi_n\to f$ in $\mu$-measure and $\Su|\psi_n-f|\,d\mu\to0$; moreover $\psi_n\ge0$ if $f\ge0$, and if $|f|\le M$, then $|\psi_n|\le M$ and $\Ch|\psi_n-f|\,d\mu\to0$;
\item the function $f\equiv1$ is a limit in $\mu$-measure of functions in $C_c(X)$;
\item $\mu(X\setminus U)=0$.
\end{enumerate}
If moreover $\mu$ is null-additive and \textup{(c)} holds, then every Borel $f\ge0$ with $\Ch f\,d\mu<\infty$ is a limit in the Choquet mean of functions in $C_c^+(X)$, i.e., there are $\psi_k\in C_c^+(X)$ with $\Ch|\psi_k-f|\,d\mu\to0$.
\end{theorem}

\begin{corollary}\label{cor:Cao-lc}
Let $X$ be a separable metric space. Theorems~4.1 and~4.2 of \cite{Cao15}, read for $f\ge0$, and with null-additivity added for unbounded Choquet integrable $f$, hold for all admissible $\mu$ if and only if $X$ is locally compact.
\end{corollary}

\begin{proof}
If $X$ is locally compact, then $U=X$ and \cref{thm:Cao-char} applies. If not, take $x_0\notin U$ and $\mu=\delta_{x_0}$, which is finite, continuous and null-additive; then $\mu(X\setminus U)=1$, so (b) of \cref{thm:Cao-char} fails for $f\equiv1\ge0$.
\end{proof}
\noindent For complete separable $X$, which is the setting of \cite{Cao15}, \cref{cor:Cao-lc} says exactly when Cao's theorems hold for all admissible measures.

\begin{proof}[Proof of \cref{thm:Cao-char}]
\emph{Preliminaries on $U$.} If $K$ is a compact neighbourhood of $x$, then every point of $\inter K$ has $K$ as a compact neighbourhood; moreover some closed ball $\bar B(x,r)\subseteq\inter K$ is compact and contained in $U$. Hence $U$ is open and locally compact. Being separable, metrizable and locally compact, $U$ is Polish \cite{Kec95}. Also $\cB(U)=\{B\cap U:B\in\cB(X)\}\subseteq\cB(X)$.\\
\emph{Step 1: compact Lusin sets in $U$.} Let $\nu:=\mu|_{\cB(U)}$, which is finite, continuous and weakly null-additive, and let $f$ be Borel and $h:=f|_U$. By \cref{lem:compact-Lusin}, applied to the Polish space $U$, the measure $\nu$ and the Borel function $h$, there are compact sets $C_1\subseteq C_2\subseteq\cdots$ in $U$ such that $h|_{C_k}$ (that is, $f|_{C_k}$) is continuous for every $k$ and $U\setminus C_k$ decreases to a set $N_U$ with $\nu(N_U)=0$.\\
\emph{Step 2: \textup{(c)}$\Rightarrow$\textup{(a)}.} We have $X\setminus C_k=(X\setminus U)\cup(U\setminus C_k)\downto(X\setminus U)\cup N_U$, which is $\mu$-null by (c) and weak null-additivity; thus $\mu(X\setminus C_k)\to0$ by continuity from above.\\
\emph{Cut-off functions.} For a compact $C\subseteq U$ there is $\chi\in C_c(X)$ with $0\le\chi\le1$ and $\chi=1$ on $C$. Indeed, every point of $C$ lies in $U$ and is therefore the centre of an open ball with compact closure; cover $C$ by finitely many such balls. Their union $V$ is open, $C\subseteq V$, and $\cl V$ is compact. If $V=X$, then $X$ is compact and we take $\chi\equiv1$. Otherwise put $\chi(x):=d(x,X\setminus V)/(d(x,X\setminus V)+d(x,C))$. The denominator is positive, since $C$ and $X\setminus V$ are disjoint closed sets; $\chi=1$ on $C$, and $\supp\chi\subseteq\cl V$ is compact. We denote by $\chi_k$ such a cut-off function for $C_k$.\\
\emph{Approximants.} Let $\tilde\varphi_k$ be a Tietze extension of $f|_{C_k}$ to $X$; replace it by $\tilde\varphi_k\vee0$ if $f\ge0$, and by $(\tilde\varphi_k\wedge M)\vee(-M)$ if $|f|\le M$. Put $\psi_k:=\chi_k\tilde\varphi_k\in C_c(X)$. Then $\psi_k=f$ on $C_k$, and the sign and the bound are preserved.\\
\emph{Estimates.} Since $\{|\psi_k-f|>0\}\subseteq X\setminus C_k$, we get $\mu(\{|\psi_k-f|\ge\alpha\})\le\mu(X\setminus C_k)$ for $\alpha>0$, $\Su|\psi_k-f|\,d\mu=\sup_{t>0}\min\bigl(t,\mu(\{|\psi_k-f|\ge t\})\bigr)\le\mu(X\setminus C_k)$, and, if $|f|\le M$, $|\psi_k-f|\le2M\1_{X\setminus C_k}$, so that $\Ch|\psi_k-f|\,d\mu\le2M\mu(X\setminus C_k)$. All these bounds tend to $0$.\\
\emph{Step 3: \textup{(b)}$\Rightarrow$\textup{(c)}.} If $\psi\in C_c(X)$ and $\psi(x)\ne0$, then $\{\psi\ne0\}$ is a neighbourhood of $x$ whose closure lies in the compact support of $\psi$; hence $x\in U$. Thus $\psi=0$ on $X\setminus U$. If $\psi_n\in C_c(X)$ converge to $1$ in $\mu$-measure, then $X\setminus U\subseteq\{|\psi_n-1|\ge\frac12\}$ and $\mu(X\setminus U)\le\mu(\{|\psi_n-1|\ge\frac12\})\to0$.\\
\textup{(a)}$\Rightarrow$\textup{(b)} is trivial.\\
\emph{Step 4: unbounded $f$ in the Choquet case, under null-additivity.} Let $f\ge0$ with $\Ch f\,d\mu<\infty$, and put $T(s):=\int_s^\infty\mu(\{f>r\})\,dr$, so that $T(s)\to0$ as $s\to\infty$ (\cref{lem:elementary}(g)). Fix $\eta>0$ and $M>0$ with $T(M)<\eta$. The set $N:=\bigcap_k(X\setminus C_k)$ is $\mu$-null by Step~2. For $t\in[0,2M]$ the functions $h_k(t):=\mu(\{f>M+t\}\cup(X\setminus C_k))$ decrease, as $k\to\infty$, to $\mu(\{f>M+t\}\cup N)=\mu(\{f>M+t\})$, by continuity from above and null-additivity. Each $h_k$ is non-increasing in $t$, hence measurable, and bounded by $\mu(X)$. Therefore,
\[
\int_0^{2M}h_k(t)\,dt\longrightarrow\int_M^{3M}\mu(\{f>s\})\,ds\le T(M),
\]
and we can choose $k$ with $\int_0^{2M}h_k<2\eta$. Let $f_M:=f\wedge M$, let $\tilde\varphi$ be a Tietze extension of $f_M|_{C_k}$, and put $\psi:=\chi_k\cdot\bigl((\tilde\varphi\wedge M)\vee0\bigr)$. Then $\psi\in C_c^+(X)$, $\psi=f_M$ on $C_k$ and $0\le\psi\le M$. On $C_k$ we have $|f-\psi|=(f-M)^+$, and off $C_k$ we have $|f-\psi|\le f+M$. Consequently
\begin{align*}
\{|f-\psi|>t\}&\subseteq\{f>M+t\}\cup(X\setminus C_k) &&(t<2M),\\
\{|f-\psi|>t\}&\subseteq\{f>t-M\} &&(t\ge2M),
\end{align*}
and
\[
\Ch|f-\psi|\,d\mu\le\int_0^{2M}h_k(t)\,dt+\int_{2M}^\infty\mu(\{f>t-M\})\,dt<2\eta+T(M)<3\eta.
\]
Since $\eta$ was arbitrary, this completes the proof.
\end{proof}
\noindent Thus, for a given $\mu$, the conclusions of \cite[Theorems~4.1--4.2]{Cao15} for convergence in measure, for the Sugeno integral and for the Choquet integral of bounded functions hold exactly when $\mu$ is concentrated on the open set $U$; for unbounded Choquet integrable $f\ge0$ we obtain them when, in addition, $\mu$ is null-additive, and the general weakly null-additive case is open (\cref{Q1}). On locally compact spaces we also obtain step-function approximation on intervals.

\begin{theorem}\label{thm:Cao-corrected}
Let $X$ be a locally compact separable metric space, $\mu$ a finite, continuous, weakly null-additive monotone measure on $\cB$, and $f$ a Borel function.
\begin{enumerate}[(a)]
\item There are $\psi_n\in C_c(X)$ with $\psi_n\to f$ in $\mu$-measure and $\Su|\psi_n-f|\,d\mu\to0$. If $f\ge0$, then $\psi_n\ge0$ (so that $\psi_n\in C_0^+(X)$ in the notation of \cite{Cao15}). If $|f|\le M$, then $|\psi_n|\le M$ and $\Ch|\psi_n-f|\,d\mu\le2M\mu(X\setminus C_n)\to0$, where $C_n$ are the sets of \cref{lem:compact-Lusin}.
\item If $\mu$ is moreover null-additive (then weak null-additivity is automatic and may be dropped from the hypotheses) and $\Ch|f|\,d\mu<\infty$, then there are $\psi_k\in C_c(X)$, with $\psi_k\ge0$ if $f\ge0$, such that $\Ch|\psi_k-f|\,d\mu\to0$.
\item \textup{(Step functions; \cite[Theorem~9]{LY04}.)} Let $X=[a,b]$, $\mu$ a finite, continuous, null-additive monotone measure on $\cB([a,b])$, and $f$ a Borel function with $\Ch|f|\,d\mu<\infty$. Then there are step functions $s_n$ (finite linear combinations of indicators of intervals) with $\Ch|s_n-f|\,d\mu\to0$. The version for the Sugeno integral and for convergence in measure holds for weakly null-additive $\mu$ and every Borel $f$.
\end{enumerate}
\end{theorem}

\begin{proof}
Since $X$ is locally compact, $U=X$, and $X$ is Polish (see the preliminaries in the proof of \cref{thm:Cao-char}), so \cref{lem:compact-Lusin} applies; let $(C_k)$ be as there. For each $k$ let $\chi_k\in C_c(X)$ be a cut-off function for $C_k$, constructed as in Step~2 of the proof of \cref{thm:Cao-char}: $0\le\chi_k\le1$ and $\chi_k=1$ on $C_k$.\\
(a) This is the implication (c)$\Rightarrow$(a) of \cref{thm:Cao-char}, since $\mu(X\setminus U)=0$; with $U=X$, the compact sets in Step~1 of its proof are those of \cref{lem:compact-Lusin} applied to $X$. Explicitly, let $\tilde\varphi_n$ be a Tietze extension of $f|_{C_n}$. For $f\ge0$ we replace it by $\tilde\varphi_n\vee0$, and for $|f|\le M$ by $(\tilde\varphi_n\wedge M)\vee(-M)$; in each case the result still equals $f$ on $C_n$. Put $\psi_n:=\chi_n\tilde\varphi_n$. Then $\psi_n\in C_c(X)$, $\psi_n=f$ on $C_n$, the sign and the bound are preserved, and $\{|\psi_n-f|>0\}\subseteq X\setminus C_n$. We conclude exactly as in \cref{thm:QC-approx}(a): $\mu(\{|\psi_n-f|\ge\alpha\})\le\mu(X\setminus C_n)\to0$, $\Su|\psi_n-f|\,d\mu\le\mu(X\setminus C_n)$, and, if $|f|\le M$, $|\psi_n-f|\le2M\1_{X\setminus C_n}$.\\
(b) We repeat the proof of \cref{thm:QC-approx}(c) with the increasing compact sets $C_k$ in place of $F'_k$. The only facts used there are:
(1) $X\setminus C_k$ decreases to a $\mu$-null set $N$ (\cref{lem:compact-Lusin});
(2) $\mu$ is finite and continuous from above, and null-additivity gives $\mu(\{|f|>M+t\}\cup N)=\mu(\{|f|>M+t\})$;
(3) the approximant $\psi$ satisfies $\psi=f_M$ on $C_k$ and $|\psi|\le M$.
For (3) we take $\psi:=\chi_k\cdot((\tilde\varphi\wedge M)\vee(-M))$, where $\tilde\varphi$ is a Tietze extension of $f_M|_{C_k}$ and $\chi_k$ is a cut-off function for $C_k$; for $f\ge0$ we use $\psi\vee0$. The inclusions for the level sets in the proof of \cref{thm:QC-approx}(b) hold verbatim (off $C_k$ only $|\psi|\le M$ was used), and we obtain $\Ch|f-\psi|\,d\mu<3\eta$.\\
(c) Let $\eta>0$. Since $X=[a,b]$ is compact, (b) gives $\varphi\in C([a,b])$ with $\nC{f-\varphi}<\eta$. (For the Sugeno version we use (a) instead, with $\varphi$ continuous and $\mu(\{|f-\varphi|\ge\alpha\})$ small.) The function $\varphi$ is uniformly continuous, so there is a step function $s$ with $\sup|\varphi-s|\le\eta':=\eta/(1+\mu(X))$. Since $|f-s|\le|f-\varphi|+\eta'$, we have $\{|f-s|>t\}\subseteq\{|f-\varphi|>t-\eta'\}$ for $t\ge\eta'$. Hence
\[
\nC{f-s}\le\int_0^{\eta'}\mu(X)\,dt+\int_{\eta'}^\infty\mu(\{|f-\varphi|>t-\eta'\})\,dt\le\eta'\mu(X)+\nC{f-\varphi}<2\eta .
\]
Likewise $\mu(\{|f-s|\ge2\alpha\})\le\mu(\{|f-\varphi|\ge\alpha\})$ once $\eta'\le\alpha$, which gives the version for convergence in measure, and the Sugeno version by \cref{lem:elementary}(i).
\end{proof}

\begin{remark}\label{rem:Cao}
(1) Theorem~9 of \cite{LY04} asserts (c) for weakly null-additive $\mu$; its Choquet part for unbounded $f$ is stated without proof (it is introduced by ``Similarly, we can obtain the following result''). \Cref{thm:Cao-corrected}(c) proves it under null-additivity. The case of a weakly null-additive $\mu$ which is not null-additive remains open (\cref{sec:concl}).\\
(2) Corollary~3 of \cite{LY04} already notes that Theorems~5--9 of \cite{LY04} hold for null-additive $\mu$. That corollary still rests on the unbounded case of Theorem~9, which is now supplied by \cref{thm:Cao-corrected}(c).\\
(3) Local compactness is necessary for approximation by functions in $C_c(X)$ for all admissible $\mu$ (\cref{cor:Cao-lc}); for a single $\mu$, what is needed is $\mu(X\setminus U)=0$ (\cref{thm:Cao-char}). Explicitly, suppose that some $x_0\in X$ has no compact neighbourhood, and take $\mu=\delta_{x_0}$ (finite and $\sigma$-additive, hence continuous and null-additive) and $f\equiv1$. Every $\psi\in C_c(X)$ satisfies $\psi(x_0)=0$; otherwise $\{\psi\ne0\}$ would be a neighbourhood of $x_0$ with compact closure, since it lies in the compact support. Hence $|\psi-1|(x_0)=1$ and $\Ch|\psi-1|\,d\delta_{x_0}=\Su|\psi-1|\,d\delta_{x_0}=\delta_{x_0}(\{|\psi-1|\ge\sigma\})=1$ for $\sigma\le1$. Thus \cite[Theorems~4.1--4.2]{Cao15} fail on every complete separable metric space which is not locally compact; \cref{ex:Cao} is the case $X=\ell^2$, $x_0=0$. Without local compactness, \cref{thm:QC-approx} and \cref{cor:unbounded} still give approximants in $C_b(X)$ on arbitrary metric spaces.
\end{remark}

%==========================================================================
\section{The Lusin property from regularity}\label{sec:lusin}

The scheme ``regularity gives the Lusin property for simple functions, and an Egoroff-type property passes it to limits'' is standard. It appears in the work of Wu and Ha \cite{WH94} for finite autocontinuous fuzzy measures and of Jiang and Suzuki \cite{JS96}, and in \cite[proof of Theorem~5.3, step (a)]{LM11}, which coincides with Step~1 below. In \cite[proof of Theorem~3 and Remark~6]{Kaw07} it is shown, for a Riesz space $V$ which is order separable and has the multiple Egoroff property, that Lusin's theorem holds for every weakly null-additive, regular, $V$-valued fuzzy Borel measure (that is, continuous from above and from below) on a normal space. In \cite[Theorem~3.1]{CWL16} the Lusin property for simple functions is obtained from regularity and weak asymptotic null-additivity. For capacities, the equivalence ``every Borel function is quasicontinuous if and only if every Borel set is quasi-open and quasi-closed'' is folklore. The only new point in the following proposition is the use of property \Sone{} (with the p.g.p.) in place of condition (E), of continuity, or of an Egoroff property of the range space.

\begin{proposition}\label{prop:regular-Lusin}
Let $\mu$ be a monotone measure on $(X,\cA)$ with the p.g.p.\ and property \Sone, and suppose that every $A\in\cA$ is regular. Then every $f\in\M$ is $\mu$-quasicontinuous. In particular, for $\cA=\cB$ and $\mu$ with the p.g.p.\ and \Sone, the Lusin property \Lu$[\mu]$ holds if and only if $\mu$ is regular.
\end{proposition}

\begin{proof}
\emph{Step 1 (finitely-valued $f$).} Let $f=\sum_{i\le m}c_i\1_{A_i}$, where $(A_i)$ is a finite partition of $X$ into sets of $\cA$ and the $c_i$ are distinct. Given $\varepsilon>0$, let $\delta>0$ be such that the union of $m$ sets of measure less than $\delta$ has measure less than $\varepsilon$ (\cref{rem:pgp}), and take $F_i\subseteq A_i\subseteq G_i$ with $F_i\in\cC$, $G_i\in\cO$ and $\mu(G_i\setminus F_i)<\delta$. The set $\Phi:=\bigcup_iF_i$ is closed, and $f$ is locally constant on $\Phi$, since the $F_i$ are finitely many disjoint closed sets, each relatively open in $\Phi$. Moreover $X\setminus\Phi=\bigcup_i(A_i\setminus F_i)\subseteq\bigcup_i(G_i\setminus F_i)$, so $\mu(X\setminus\Phi)<\varepsilon$, and $G:=X\setminus\Phi$ is open.\\
\emph{Step 2 (bounded $f$).} The functions $s_k:=2^{-k}\lfloor2^kf\rfloor$ are finitely-valued, $\cA$-measurable, and $|f-s_k|\le2^{-k}$. By Step~1 there are open sets $O_k$ with $\mu(O_k)<1/k$ such that $s_k|_{X\setminus O_k}$ is continuous. By \Sone{} there are $k_i$ with $\mu(U_j)\to0$, where $U_j:=\bigcup_{i\ge j}O_{k_i}$ is open. On $X\setminus U_j$ every $s_{k_i}$ with $i\ge j$ is continuous, and $s_{k_i}\to f$ uniformly; so $f|_{X\setminus U_j}$ is continuous, and $\mu(U_j)$ is as small as we like.\\
\emph{Step 3 (general $f$).} The function $g:=\arctan\circ f$ is bounded and $\cA$-measurable, hence quasicontinuous by Step~2; on the same closed sets, $f=\tan\circ g$ is continuous.\\
The last statement follows from \cref{lem:indicator-QC}, which shows that \Lu$[\mu]$ implies regularity.
\end{proof}

It seems natural to ask whether, on a compact metric space, every regular and order continuous finite monotone measure has the Lusin property. \Cref{prop:regular-Lusin} answers the analogous question with \Sone{} and the p.g.p.\ in place of order continuity and compactness. For instance, for a countably subadditive $\mu$ (such as $\Hd$ on $2^{\R^n}$), every $f\in\M$ is $\mu$-quasicontinuous if and only if $\mu$ is regular.

%==========================================================================
\section{Density of continuous functions}\label{sec:density}

\begin{theorem}[Indicator functions]\label{thm:indicator}
Let $\mu$ be an arbitrary monotone measure on $(X,\cA)$ and $A\in\cA$. The following are equivalent:
\begin{enumerate}[(i)]
\item $\1_A\in\clM(C(X))$;
\item $\1_A\in\clS(C(X))$;
\item $\1_A\in\clC(C_b(X))$, with approximants $0\le\varphi_n\le1$;
\item $\1_A\in\clC(C(X))$;
\item $\1_A\in\clSh(C(X))$;
\item $A$ has property \RD.
\end{enumerate}
Moreover, in (iii) one can choose $\varphi$ with $|\1_A-\varphi|\le\1_E$, where $E$ is the set in \RD. Every regular set satisfies (vi), and, under \OS, (vi) implies that $A$ is regular (\cref{lem:implications}(g)).
\end{theorem}

\begin{proof}
(vi)$\Rightarrow$(iii): Given $\varepsilon>0$, take $F\subseteq G$ and $E$ as in \RD. Define $\varphi\equiv0$ if $F=\emptyset$, $\varphi\equiv1$ if $G=X$, and otherwise
\[
\varphi(x):=\frac{d(x,X\setminus G)}{d(x,X\setminus G)+d(x,F)} .
\]
The denominator vanishes only on $F\cap(X\setminus G)=\emptyset$, so $\varphi$ is continuous (even locally Lipschitz), $0\le\varphi\le1$, $\varphi=1$ on $F$ and $\varphi=0$ on $X\setminus G$. Let $x\notin E$. If $x\in A$, then $x\notin A\setminus F$, so $x\in F$ and $\varphi(x)=1$; if $x\notin A$, then $x\notin G\setminus A$, so $x\notin G$ and $\varphi(x)=0$. Hence $|\1_A-\varphi|\le\1_E$, and by monotonicity $\Ch|\1_A-\varphi|\,d\mu\le\mu(E)<\varepsilon$, $\Sh|\1_A-\varphi|\,d\mu\le\mu(E)$, $\Su|\1_A-\varphi|\,d\mu\le\min(1,\mu(E))$ and $\mu(\{|\1_A-\varphi|\ge\alpha\})\le\mu(E)$.\\
(iii)$\Rightarrow$(iv)$\Rightarrow$(i), (v)$\Rightarrow$(i) and (ii)$\Leftrightarrow$(i) follow from \cref{lem:elementary}(i), and (iii)$\Rightarrow$(v) holds because the Shilkret integral is dominated by the Choquet integral.\\
(i)$\Rightarrow$(vi): Given $\varepsilon>0$, take $\varphi\in C(X)$ with $\mu(E')<\varepsilon$, where $E':=\{|\1_A-\varphi|\ge1/3\}\in\cA$. Put $G:=\{\varphi>1/3\}\in\cO$ and $F:=\{\varphi\ge2/3\}\in\cC$; then $F\subseteq G$. If $x\in A\setminus G$, then $|1-\varphi(x)|\ge2/3$; if $x\in G\setminus A$, then $|\varphi(x)|>1/3$; if $x\in A\setminus F$, then $|1-\varphi(x)|>1/3$; and if $x\in F\setminus A$, then $|\varphi(x)|\ge2/3$. Hence $(A\triangle F)\cup(A\triangle G)\subseteq E'$.
\end{proof}
\noindent No finiteness of $\mu(A)$, no subadditivity and no continuity of $\mu$ is used here. If $\mu(A)<\infty$, then $\1_A\in\Lone(\mu)$; whether the approximants lie in $\Lone(\mu)$ depends on $\mu$ (under \SK{} they do, since $\{\varphi>t\}\subseteq A\cup E$).

\begin{theorem}[Density in the Choquet space]\label{thm:Choquet-density}
Let $\mu$ be a monotone measure on $(X,\cA)$ and consider the statements:
\begin{enumerate}[(a)]
\item $C_b(X)\cap\Lone(\mu)$ is dense in $(\Lone(\mu),\nC\cdot)$, that is, every $f\in\Lone(\mu)$ is a $\nC\cdot$-limit of such functions;
\item every $f\in\Lone(\mu)$ lies in $\clC(C(X))$;
\item every $A\in\cA$ with $\mu(A)<\infty$ has property \RD.
\end{enumerate}
Then (a)$\Rightarrow$(b)$\Rightarrow$(c) for every $\mu$, and (c)$\Rightarrow$(a) if $\mu$ satisfies \SK{} (e.g.\ if $\mu$ is a submeasure, satisfies \UA, or is $K$-relaxed subadditive).
\end{theorem}

\begin{proof}
(a)$\Rightarrow$(b) is trivial. (b)$\Rightarrow$(c): if $\mu(A)<\infty$, then $\1_A\in\Lone(\mu)$ by \cref{lem:elementary}(e), and \cref{thm:indicator} (iv)$\Rightarrow$(vi) applies.\\
(c)$\Rightarrow$(a). \emph{Step 1 (for every $\mu$): simple functions with finitely many values and level sets of finite measure are dense.} For $f\in\Lone(\mu)$ and $n\in\N$ put $s_n:=\sgn(f)\cdot2^{-n}\lfloor2^n\min(|f|,n)\rfloor$. Then $s_n$ is $\cA$-measurable, takes finitely many values, $|s_n|\le|f|$ and $\sgn s_n\in\{0,\sgn f\}$, so $|f-s_n|=|f|-|s_n|$. On $\{|f|\le n\}$ we have $0\le|f|-|s_n|<2^{-n}$, and on $\{|f|>n\}$ we have $|s_n|=n$. Hence $\{|f-s_n|>t\}\subseteq\{|f|>t\}$ for $t<2^{-n}$ and $\{|f-s_n|>t\}\subseteq\{|f|>n+t\}$ for $t\ge2^{-n}$, and
\[
\nC{f-s_n}\le\int_0^{2^{-n}}\mu(\{|f|>t\})\,dt+\int_n^\infty\mu(\{|f|>s\})\,ds\longrightarrow0,
\]
since $t\mapsto\mu(\{|f|>t\})$ is Lebesgue integrable on $(0,\infty)$. Each nonzero level set of $s_n$ lies in $\{|f|\ge2^{-n}\}$, whose measure is at most $2^n\nC f<\infty$ (\cref{lem:elementary}(b)).\\
\emph{Step 2.} Fix $\eta>0$ and $s:=s_n$ with $\nC{f-s}<\eta$, and write $s=\sum_{i=1}^mc_i\1_{A_i}$ with disjoint $A_i$, $c_i\ne0$ and $\mu(A_i)<\infty$; put $c:=\sum_i|c_i|$. Let $K$ and $\delta(\cdot)$ be as in \SK, $\varepsilon:=\eta/(2c)$, and let $\delta'>0$ be such that the union of $m$ sets of measure less than $\delta'$ has measure less than $\delta(\varepsilon)$ (\SK{} implies the p.g.p.\ by \cref{lem:implications}(c); see \cref{rem:pgp}). By (c) and \cref{thm:indicator} there are $\varphi_i\in C_b(X)$ with $0\le\varphi_i\le1$ and $E_i\in\cA$ with $\mu(E_i)<\delta'$ and $|\1_{A_i}-\varphi_i|\le\1_{E_i}$. Put $\varphi:=\sum_ic_i\varphi_i\in C_b(X)$ and $E:=\bigcup_iE_i$, so that $\mu(E)<\delta(\varepsilon)$. Off $E$ we have $\varphi=s$, and $|s-\varphi|\le c$ everywhere.\\
\emph{Step 3.} We have $\{|f-\varphi|>t\}\subseteq\{|f-s|>t\}\cup E$ for all $t$, and $\{|f-\varphi|>t\}\subseteq\{|f-s|>t-c\}\subseteq\{|f-s|>t/2\}$ for $t\ge2c$. Hence
\begin{align*}
\nC{f-\varphi}&\le\int_0^{2c}\bigl[K\mu(\{|f-s|>t\})+\varepsilon\bigr]dt+\int_{2c}^\infty\mu(\{|f-s|>t/2\})\,dt\\
&\le K\nC{f-s}+2c\varepsilon+2\nC{f-s}<(K+3)\eta .
\end{align*}
Finally $\varphi\in\Lone(\mu)$: $\{|\varphi|>t\}\subseteq\{|s|>t\}\cup E$, and this set is empty for $t\ge c$, so $\nC\varphi\le K\nC s+c\varepsilon\le K\nC f+c\varepsilon<\infty$.
\end{proof}

\begin{remark}\label{rem:simple}
Kawabe and Yamada \cite{KY23} study the space of measurable functions and the Choquet--Lorentz spaces; the space $\Lone(\mu)$ is the Choquet--Lorentz space with exponents $(1,1)$. Step~1 above shows that, for $\Lone(\mu)$ with the prenorm $\nC\cdot$, simple functions with finitely many values and level sets of finite measure are dense for every monotone measure, without any continuity assumption. The reason is the tail integrability $\int_n^\infty\mu(\{|f|>s\})\,ds\to0$ of Choquet integrable functions. For the space defined by the Shilkret integral (of weak type) no such statement can hold, even for the Lebesgue measure (\cref{prop:shilkret}).
\end{remark}

\begin{theorem}[Density in measure and in the Sugeno integral]\label{thm:measure-density}
Let $\mu$ be a monotone measure on $(X,\cA)$. Then
\begin{enumerate}[(a)]
\item If $\clM(C_b(X))=\M$ (equivalently, $\clS(C_b(X))=\M$), then $\mu$ is order continuous and every $A\in\cA$ has property \RD.
\item If $\mu$ has the p.g.p., is order continuous, and every $A\in\cA$ has property \RD, then $\clM(C_b(X))=\clS(C_b(X))=\M$.
\item Order continuity is not necessary for $\clM(C(X))=\M$ (with unbounded approximants): let $X=\N$ with the discrete metric and $\mu(A)=1$ if $A$ is infinite, $\mu(A)=0$ otherwise. Every function is continuous, $\mu$ is a submeasure (so it has the p.g.p.), $\mu$ is not order continuous, and $C_b$ is not dense in measure: for $f(n)=n$ and $|\varphi|\le M$, the set $\{|f-\varphi|\ge1\}\supseteq\{n\ge M+1\}$ is infinite.
\end{enumerate}
\end{theorem}

\begin{proof}
(a) Property \RD{} follows from \cref{thm:indicator}. For order continuity, let $A_n\downto\emptyset$ and $f:=\sum_n\1_{A_n}$. For each $x$, $N(x):=\min\{n:x\notin A_n\}$ exists and $x\in A_n$ if and only if $n<N(x)$; so $f(x)=N(x)-1<\infty$, $\{f\ge k\}=A_k$ for $k\ge1$, and $f\in\M$. Given $\varepsilon>0$, pick $\varphi\in C_b(X)$, $|\varphi|\le M$, with $\mu(\{|f-\varphi|\ge1\})<\varepsilon$. For $n\ge M+1$ and $x\in A_n$ we have $f(x)\ge n\ge M+1$, hence $|f(x)-\varphi(x)|\ge1$. Thus $A_n\subseteq\{|f-\varphi|\ge1\}$ and $\mu(A_n)<\varepsilon$ for $n\ge M+1$.\\
(b) Let $f\in\M$ and $\alpha,\varepsilon>0$. We make the choices in the following order. (1) Let $\delta_2:=\delta_{\rm pgp}(\varepsilon)$, so that two sets of measure less than $\delta_2$ have union of measure less than $\varepsilon$. (2) Since $f$ is finite-valued, $\{|f|>N\}\downto\emptyset$ as $N\to\infty$; by order continuity fix $N$ with $\mu(\{|f|>N\})<\delta_2$. (3) Put $f_N:=(f\wedge N)\vee(-N)$ and $s:=(\alpha/2)\lfloor2f_N/\alpha\rfloor$, a finitely-valued $\cA$-measurable function with $|f_N-s|\le\alpha/2<\alpha$; write $s=\sum_{i\le m}c_i\1_{A_i}$ with disjoint $A_i$ and $c_i\ne0$. Now $m$ is fixed. (4) Choose $\delta'>0$ such that the union of $m$ sets of measure less than $\delta'$ has measure less than $\delta_2$ (\cref{rem:pgp}), and, by \RD{} and \cref{thm:indicator}, $\varphi_i\in C_b(X)$ with $0\le\varphi_i\le1$ and $E_i\in\cA$ with $\mu(E_i)<\delta'$ and $|\1_{A_i}-\varphi_i|\le\1_{E_i}$. Put $\varphi:=\sum_ic_i\varphi_i\in C_b(X)$ and $E:=\bigcup_iE_i$, so that $\mu(E)<\delta_2$. Off $E$ we have $\varphi=s$, and off $\{|f|>N\}$ we have $f=f_N$. So $|f-\varphi|=|f_N-s|<\alpha$ on $X\setminus(\{|f|>N\}\cup E)$, that is, $\{|f-\varphi|\ge\alpha\}\subseteq\{|f|>N\}\cup E$, whose measure is less than $\varepsilon$ by (1). Taking $\alpha=\varepsilon=1/k$ we obtain $\varphi_k\in C_b(X)$ with $\varphi_k\to f$ in measure; the Sugeno case follows from \cref{lem:elementary}(i).
\end{proof}

Finiteness of $\mu$ is not needed in (b). The proof of (a) uses only $|\varphi|\le M$ and not the continuity of $\varphi$. It therefore shows: \emph{if the bounded $\cA$-measurable functions are dense in $\M$ with respect to convergence in $\mu$-measure, then $\mu$ is order continuous.} This seems to complement the density results for simple functions in the space of measurable functions studied in \cite{KY23}, where continuity assumptions of this kind play a role; we refer the reader to \cite{KY23} for the precise statements. If \Lu$[\mu]$ holds, then $\clM(C(X))=M(\cB)$ for every $\mu$ (\cref{thm:QC-approx}(a)). A characterization of $\clM(C(X))=\M$ without order continuity is open (\cref{sec:concl}).

\begin{corollary}[Density and the Lusin property]\label{cor:density-Lusin}
Let $\mu$ be a monotone measure on $(X,\cA)$ with \Sone, \OS{} and \SK. The following are equivalent:
\begin{enumerate}[(a)]
\item $C_b(X)\cap\Lone(\mu)$ is dense in $\Lone(\mu)$;
\item every $f\in\Lone(\mu)$ is $\mu$-quasicontinuous;
\item every $A\in\cA$ with $\mu(A)<\infty$ is regular;
\item every $A\in\cA$ with $\mu(A)<\infty$ has property \RD.
\end{enumerate}
If $\mu$ is finite, these are also equivalent to: (e) every $f\in\M$ is $\mu$-quasicontinuous (that is, \Lu$[\mu]$ for $\cA=\cB$); and (f) $\mu$ is regular.
\end{corollary}

\begin{proof}
(a)$\Rightarrow$(b): \cref{thm:approx-QC}. (b)$\Rightarrow$(c): $\1_A\in\Lone(\mu)$ and \cref{lem:indicator-QC}. (c)$\Rightarrow$(d): \cref{lem:implications}(g). (d)$\Rightarrow$(a): \cref{thm:Choquet-density}. If $\mu$ is finite, then every bounded $f\in\M$ belongs to $\Lone(\mu)$, so (b) gives quasicontinuity of bounded functions, and (e) follows as in Step~3 of the proof of \cref{prop:regular-Lusin}; (e)$\Rightarrow$(b) is trivial, and (f)$\Leftrightarrow$(c) because all sets have finite measure.
\end{proof}

We note a consequence for submeasures. Let $\mu$ be a finite order continuous submeasure on a metric space $X$. If $A_n\downarrow A$, then $\mu(A)\le\mu(A_n)\le\mu(A)+\mu(A_n\setminus A)$ and $A_n\setminus A\downarrow\emptyset$; if $A_n\uparrow A$, then $\mu(A_n)\le\mu(A)\le\mu(A_n)+\mu(A\setminus A_n)$ and $A\setminus A_n\downarrow\emptyset$. Hence $\mu$ is continuous from above and from below, and consequently countably subadditive, since $\mu(\bigcup_{i}A_i)=\lim_n\mu(\bigcup_{i\le n}A_i)\le\sum_i\mu(A_i)$. In particular $\mu$ is a finite weakly null-additive fuzzy measure, so \Lu$[\mu]$ holds by \cite[Theorem~4]{LY04}. Moreover $\mu$ has \SK{} and \Sone, and it is regular (by \Lu$[\mu]$ and \cref{lem:indicator-QC}), so it has \OS{} (\cref{lem:implications}(f)). Hence $C_b(X)$ is dense in $\Lone(\mu)$ for every order continuous finite submeasure on a metric space, in particular on a compact one.

\begin{proposition}[No density for the Shilkret integral]\label{prop:shilkret}
Let $X=[0,1]$, $\mu=m$ the Lebesgue measure, and $f(x)=1/x$ for $x>0$, $f(0)=0$. Then $\Sh|f|\,dm=1$ and $f$ is $m$-quasicontinuous, but $\Sh|f-\varphi|\,dm\ge1$ for every $\varphi\in C([0,1])$.
\end{proposition}

This is the classical fact that bounded functions are not dense in weak $L^1$ (see e.g.\ \cite[Section~1.4]{Gra14}); we include it only to explain why the Shilkret integral is excluded from the density statements.

\begin{proof}
$\Sh|f|\,dm=\sup_tt\,m(\{1/x\ge t\})=\sup_tt\min(1,1/t)=1$, and $f$ is continuous on $[\eta,1]$ with $m([0,\eta))=\eta$. Let $\varphi\in C([0,1])$, so $|\varphi|\le M$ for some $M\ge0$. For $t>0$ with $M+t\ge1$ we have $\{|f-\varphi|\ge t\}\supseteq(0,1/(M+t)]$, hence $\Sh|f-\varphi|\,dm\ge\sup_{t\ge1}t/(M+t)=1$.
\end{proof}

Thus quasicontinuity does not imply approximability in the Shilkret integral, and $C([0,1])$ is not dense in the Shilkret space of $m$; bounded quasicontinuous functions are approximable by \cref{thm:QC-approx}(a). The Shilkret integral therefore appears only in \cref{prop:riesz,prop:riesz-au,thm:approx-QC,thm:indicator}.

%==========================================================================
\section{Hausdorff content}\label{sec:hausdorff}

Following \cite[Section~2]{HH25}, for $0<\delta\le n$ and $E\subseteq\R^n$ we put
\[
\Hd(E):=\inf\Bigl\{\sum_i\ell(Q_i)^\delta:\ E\subseteq\inter\Bigl(\bigcup_iQ_i\Bigr)\Bigr\},
\]
where the infimum is taken over finite or countable families of half-open dyadic cubes $Q_i$ and $\ell(Q)$ denotes the side length of $Q$. We also use $\Hstd(E):=\inf\{\sum_i\diam(U_i)^\delta:\ E\subseteq\bigcup_iU_i\}$, with arbitrary sets $U_i$. Only the following elementary comparison is needed:
\begin{equation}\label{eq:compare}
\Hstd(E)\le n^{\delta/2}\,\Hd(E),
\end{equation}
since an admissible dyadic family covers $E$ and $\diam Q=\sqrt n\,\ell(Q)$.

\begin{theorem}\label{thm:cantor}
Let $n\ge1$ and $0<\delta<n$. Put $s:=2^{-n/\delta}\in(0,1/2)$, let $C_s\subseteq[0,1]$ be the symmetric Cantor set (at step $k$ there are $2^k$ closed intervals of length $s^k$, and each interval $J$ of step $k$ is replaced by its two end subintervals of length $s^{k+1}$), and let $K:=C_s^n\subseteq[0,1]^n$. Then
\begin{enumerate}[(a)]
\item $K$ is compact and Lebesgue-null, and $0<\Hd(K)<\infty$; so $\1_K\in\Lone(\R^n,\Hd)$ and $\1_K=0$ Lebesgue-almost everywhere;
\item $K$ is not $\Hd$-regular; hence it does not have property \RD, and $\1_K\notin\clM(C(\R^n))$. In particular, $\1_K$ is not an $\mathrm{N}L^1$-limit of continuous functions.
\end{enumerate}
Consequently, $C_b(\R^n)$ and $C(\R^n)$ are not dense in $\Lone(\R^n,\cA,\Hd)$ for $\cA=\cB$, for the Lebesgue $\sigma$-algebra, and for $\cA=2^{\R^n}$.
\end{theorem}
\noindent As explained in the introduction, the non-density statement itself is stated in \cite[Remark~6.3]{HH25}; the examples given there cover $0<\delta\le n-1$ (a ball) and $n=2$, $1<\delta<2$ (a modified von Koch snowflake). The example of \cref{thm:cantor} covers all cases at once, in particular $n=1$ and $n\ge3$ with $n-1<\delta<n$.
\begin{proof}
A \emph{level-$k$ cell} is a product $Q=J_1\times\dots\times J_n$ of step-$k$ intervals; there are $2^{nk}$ of them, each of side $s^k$, and $K=\bigcap_k\bigcup\{\text{level-}k\text{ cells}\}$. Let $\nu$ be the natural probability measure on $K$ (the product of the Cantor measures), which gives mass $2^{-nk}$ to each level-$k$ cell. Note that $s^\delta=2^{-n}$, i.e.\ $2^{-nk}=(s^k)^\delta$.\\
\emph{Fact A (projections).} An interval of length less than $s^j$ meets at most two step-$j$ intervals, since step-$j$ intervals are disjoint and of length $s^j$, and an interval meeting three of them would contain the middle one. Hence a set $B\subseteq\R^n$ of diameter $d<s^j$ meets at most $2^n$ level-$j$ cells, because its coordinate projections have length at most $d$.\\
\emph{Fact B (growth of $\nu$).} $\nu(B)\le4^n\diam(B)^\delta$ for every $B$. If $d:=\diam B\ge1$, then $\nu(B)\le1\le d^\delta$. If $0<d<1$, choose $j\ge0$ with $s^{j+1}\le d<s^j$; by Fact~A, $\nu(B)\le2^n2^{-nj}=2^n(s^j)^\delta\le2^n(d/s)^\delta=4^nd^\delta$. If $d=0$, then $\nu(B)\le2^n2^{-nj}$ for all $j$, so $\nu(B)=0$.\\
(a) $K$ is an intersection of compact sets, and $m(C_s)\le(2s)^k\to0$. Next, $\Hd(K)\le2^n\cdot2^\delta$, since $K\subseteq(-2,2)^n$, which is the interior of the union of the $2^n$ dyadic cubes of side $2$ having a vertex at the origin. Finally, by \eqref{eq:compare} and Fact~B, $\Hd(K)\ge n^{-\delta/2}\Hstd(K)\ge n^{-\delta/2}4^{-n}>0$: if $K\subseteq\bigcup_iU_i$, then $1=\nu(K)\le\sum_i\nu(\cl U_i)\le4^n\sum_i\diam(U_i)^\delta$.\\
(b) \emph{Holes.} For a level-$k$ cell $Q=J_1\times\dots\times J_n$ let $g_Q$ be the open middle gap of $J_1$ (of length $(1-2s)s^k$), and $H_Q:=g_Q\times J_2\times\dots\times J_n$. Points of $H_Q$ have first coordinate outside $C_s$, so $H_Q\cap K=\emptyset$; since $Q\cap K\ne\emptyset$, we have $\dist(x,K)\le\diam Q=\sqrt n\,s^k$ for $x\in H_Q$. Also $m(H_Q)=(1-2s)s^{nk}$.\\
\emph{Hole measures.} For $k\ge0$ let
\[
\sigma_k:=\sum_Q2^{-nk}\,\frac{m|_{H_Q}}{m(H_Q)}
\]
(sum over the $2^{nk}$ level-$k$ cells), a probability measure carried by $\bigcup_QH_Q$.\\
\emph{Claim C.} $\sigma_k(B)\le C_0\diam(B)^\delta$ for every Borel set $B$, where $C_0:=\max\bigl(4^n(1+2\sqrt n)^\delta,\ 2^n/(1-2s)\bigr)$ does not depend on $k$. Let $\rho:=\diam B$.
\begin{itemize}[leftmargin=2em]
\item If $\rho\ge s^k$, then $\sigma_k(B)\le\sum_{Q:\,Q\cap B\ne\emptyset}2^{-nk}=\nu\bigl(\bigcup_{Q\cap B\ne\emptyset}Q\bigr)\le\nu(B^+)$, where $B^+:=\{x:\dist(x,B)\le\sqrt n\,s^k\}$ (note that $H_Q\subseteq Q$ and $\nu(Q)=2^{-nk}$). Since $\diam B^+\le\rho+2\sqrt n\,s^k\le(1+2\sqrt n)\rho$, Fact~B gives the claim.
\item If $\rho<s^k$, then by Fact~A the set $B$ meets at most $2^n$ level-$k$ cells, hence at most $2^n$ holes; $B$ lies in a cube of side $\rho$, so $m(B\cap H_Q)\le\rho^n$. Thus
\[
\sigma_k(B)\le2^n\cdot2^{-nk}\frac{\rho^n}{(1-2s)s^{nk}}=\frac{2^n}{1-2s}\,\rho^\delta\Bigl(\frac\rho{s^k}\Bigr)^{n-\delta}\le\frac{2^n}{1-2s}\,\rho^\delta,
\]
by using $2^{-nk}=s^{k\delta}$ and $n>\delta$.
\end{itemize}
\noindent \emph{Lower bound.} Let $r>0$ and $K_r:=\{x:\dist(x,K)<r\}$. Choose $k$ with $\sqrt n\,s^k<r$. Then $\bigcup_QH_Q\subseteq K_r\setminus K$. If $K_r\setminus K\subseteq\bigcup_iU_i$, then $1=\sigma_k(\bigcup_QH_Q)\le\sum_i\sigma_k(\cl U_i)\le C_0\sum_i\diam(U_i)^\delta$. Hence $\Hstd(K_r\setminus K)\ge1/C_0$, and by \eqref{eq:compare} $\Hd(K_r\setminus K)\ge c:=n^{-\delta/2}/C_0$ for every $r>0$.\\
\emph{Non-regularity.} Let $F\subseteq K\subseteq G$ with $F$ closed and $G$ open. Since $K$ is compact, there is $r>0$ with $K_r\subseteq G$ (the function $x\mapsto\dist(x,\R^n\setminus G)$ is continuous and positive on $K$; if $G=\R^n$ any $r$ works). Then $G\setminus F\supseteq G\setminus K\supseteq K_r\setminus K$, so $\Hd(G\setminus F)\ge c$. Thus $K$ is not regular. Since $\Hd$ has \OS{} (proof of \cref{cor:HH}), $K$ does not have property \RD{} (\cref{lem:implications}(g)), and $\1_K\notin\clM(C(\R^n))$ by \cref{thm:indicator}. This holds for every $\sigma$-algebra $\cA\ni K$ to which we restrict $\Hd$, because \OS{} and regularity involve only open and closed sets and $K$. The last statement follows from (a), (b) and \cref{lem:elementary}(i).
\end{proof}

\begin{remark}\label{rem:flat}
For $0<\delta\le n-1$ a flat piece such as $[0,1]^{n-1}\times\{0\}$ already works, by the projection argument of \cref{prop:ball} below: the set $K_r\setminus K$ contains $[0,1]^{n-1}\times\{r/2\}$, whose content is bounded below independently of $r$. For $n-1<\delta<n$, on the other hand, a ball is regular, because the annuli around its boundary sphere have small $\delta$-content (an annulus of thickness $t$ has content of order $t^{\delta-n+1}$). So a fractal set is needed exactly for $n-1<\delta<n$, which includes the whole case $n=1$.
\end{remark}

\begin{proposition}[The ball; {\cite[Example~6.1]{HH25}}]\label{prop:ball}
Let $n\ge2$, $0<\delta\le n-1$, and let $B$ be the open unit ball of $\R^n$. Then $B$ is not $\Hd$-regular, and $\1_B\notin\clM(C(\R^n))$.
\end{proposition}
\noindent The result is \cite[Example~6.1]{HH25}; the following proof avoids maximal functions.

\begin{proof}
Let $F\subseteq B\subseteq G$ with $F$ closed and $G$ open. Then $F$ is compact, $F\subseteq\bar B(0,\rho_0)$ for some $\rho_0<1$, and $G\setminus F\supseteq S_\rho:=\{|x|=\rho\}$ for $\rho\in(\max(\rho_0,1/2),1)$. Let $\pi\colon\R^n\to\R^{n-1}$ drop the last coordinate; $\pi$ is $1$-Lipschitz and $\pi(S_\rho)\supseteq D$, the closed $(n-1)$-dimensional ball of radius $1/2$. Covers of $S_\rho$ are mapped to covers of $D$ with no larger diameters, so $\Hstd(S_\rho)\ge\Hstd(D)$. If $D\subseteq\bigcup_iU_i$ and $\diam U_i\le1$ for all $i$, then $\mathcal L^{n-1}(D)\le\sum_i\omega_{n-1}\diam(U_i)^{n-1}\le\omega_{n-1}\sum_i\diam(U_i)^\delta$, where $\omega_{n-1}$ is the volume of the unit ball of $\R^{n-1}$ (each $U_i$ lies in a ball of radius $\diam U_i$); otherwise $\sum_i\diam(U_i)^\delta\ge1$. Hence $\Hstd(D)\ge\min(1,2^{1-n})$ and, by \eqref{eq:compare}, $\Hd(G\setminus F)\ge n^{-\delta/2}\min(1,2^{1-n})>0$. So $B$ is not regular, and $\1_B\notin\clM(C(\R^n))$ as in the proof of \cref{thm:cantor}.
\end{proof}

\begin{corollary}[{$\delta=n$; \cite[Theorem~5.9(1)]{HH25}}]\label{cor:delta-n}
Let $\cL$ be the Lebesgue $\sigma$-algebra of $\R^n$ and $\mu:=\Hn|_{\cL}$. Then $C_b(\R^n)\cap\Lone(\mu)$ is dense in $\Lone(\R^n,\cL,\Hn)$.
\end{corollary}

\begin{proof}
(i) $m(A)\le\Hn(A)$ for $A\in\cL$, since $A\subseteq\bigcup_iQ_i$ gives $m(A)\le\sum_i\ell(Q_i)^n$. (ii) $\Hn(A)\le m(O)$ for every open $O\supseteq A$: $O$ is a countable disjoint union of half-open dyadic cubes $Q_i$, $\inter(\bigcup_iQ_i)=O\supseteq A$, and $\sum_i\ell(Q_i)^n=m(O)$. So if $\Hn(A)<\infty$, then $m(A)<\infty$, and the regularity of the Lebesgue measure gives $F\subseteq A\subseteq G$ with $m(G\setminus F)<\varepsilon$; since $G\setminus F$ is open, $\Hn(G\setminus F)\le m(G\setminus F)<\varepsilon$. Thus every $A\in\cL$ with $\mu(A)<\infty$ is regular, hence has \RD{} (\cref{lem:implications}(g)). The content is subadditive, so it satisfies \SK{} (\cref{lem:implications}(c)), and \cref{thm:Choquet-density} gives the density.
\end{proof}
\noindent This recovers \cite[Theorem~5.9(1)]{HH25} for finite-valued functions; there, approximants in $C_0$ are obtained from classical $L^1$-density. With $\cA=2^{\R^n}$, a non-measurable set $A$ with $\Hn(A)<\infty$ gives $\1_A\in\Lone$ outside $\clC(C(\R^n))$ \cite[Theorem~5.9(2)]{HH25}; by \cref{thm:indicator} this means that such a set does not have property \RD, which is consistent with our results.

%==========================================================================
\section{Further examples}\label{sec:examples}

\begin{example}[Possibility measures]\label{ex:possibility}
Let $\pi\colon X\to[0,1]$ be upper semicontinuous, $\Pi(A):=\sup_{x\in A}\pi(x)$ (with $\sup\emptyset=0$), and let $\cA\supseteq\cB$ be any $\sigma$-algebra (e.g.\ $2^X$). Then 
\begin{itemize}[leftmargin=1em]
\item $\Pi$ is maxitive: $\Pi(\bigcup_iA_i)=\sup_i\Pi(A_i)\le\sum_i\Pi(A_i)$. Hence it is countably subadditive, and has \Sone, \PS, the p.g.p., \UA{} and \SK{} with $K=1$ (\cref{lem:implications}).
\item $\Pi$ satisfies the outer regularity condition of \cref{lem:implications}(e), hence \OS: for $\eta>0$ the set $G:=\{\pi<\Pi(E)+\eta\}$ is open ($\pi$ is upper semicontinuous), $G\supseteq E$ and $\Pi(G)\le\Pi(E)+\eta$.
\item $\Ch|f|\,d\Pi=\int_0^\infty\sup\{\pi(x):|f(x)|>t\}\,dt$ and $\Su|f|\,d\Pi=\sup_x\min(|f(x)|,\pi(x))$ (interchange the two suprema:
\[\sup_t\sup_{|f(x)|\ge t}\min(t,\pi(x))=\sup_x\sup_{t\le|f(x)|}\min(t,\pi(x))\text{).}\]
\item \emph{$f$ is $\Pi$-quasicontinuous if and only if $f|_{P_\varepsilon}$ is continuous for every $\varepsilon>0$, where $P_\varepsilon:=\{\pi\ge\varepsilon\}$} (a closed set). Indeed, $\Pi(G)<\varepsilon$ forces $G\subseteq\{\pi<\varepsilon\}$, i.e.\ $P_\varepsilon\subseteq X\setminus G$. Conversely, $G:=\{\pi<\varepsilon/2\}$ is open, $\Pi(G)\le\varepsilon/2<\varepsilon$ and $X\setminus G=P_{\varepsilon/2}$.
\item Hence, by \cref{thm:closure}, $\clM(C(X))=\{f\in\M:\ f|_{P_\varepsilon}\text{ is continuous for all }\varepsilon>0\}$, and for $f\in\Lone(\Pi)$ we have $f\in\clC(C_b(X))$ if and only if $f|_{P_\varepsilon}$ is continuous for all $\varepsilon>0$.
\item \emph{Density criterion.} $C_b(X)$ is dense in $\Lone(\Pi)$ if and only if every $A\in\cA$ is regular (\cref{cor:density-Lusin}; all sets have finite measure), if and only if $\1_A|_{P_\varepsilon}$ is continuous for all $\varepsilon>0$ and all $A\in\cA$ (\cref{lem:indicator-QC}), if and only if $A\cap P_\varepsilon$ is relatively clopen in $P_\varepsilon$ for all such $A$ and $\varepsilon$, if and only if \emph{every $P_\varepsilon$ $(\varepsilon>0)$ is a discrete subspace}. (For ``only if'' note that singletons are Borel; for ``if'', every subset of a discrete space is clopen.) By \cref{prop:regular-Lusin}, this is also equivalent to \Lu$[\Pi]$ when $\cA=\cB$.
\item For instance: $\pi=\1_{\{x_0\}}$ gives density (here $\nC f=|f(x_0)|$); $\pi(x)=x$ on $[0,1]$ gives non-density ($P_\varepsilon=[\varepsilon,1]$); and $\pi\equiv1$ gives \cref{ex:nonempty}.
\end{itemize}
\end{example}

\begin{example}[$w(A)=1$ if and only if $A\ne\emptyset$]\label{ex:nonempty}
Let $X=[0,1]$ and $w(A)=1$ if $A\ne\emptyset$, $w(\emptyset)=0$; this is $\Pi$ with $\pi\equiv1$, and the conjugate of the measure $v$ of \cref{rem:S1-not-nec}. Then $\nC f=\sup|f|$ and $\Su|f|\,dw=\min(1,\sup|f|)$, so $\Lone(w)$ is the space of bounded measurable functions with the supremum norm. A function is $w$-quasicontinuous if and only if it is continuous ($w(G)<1$ forces $G=\emptyset$), and $\clC(C_b)=C[0,1]$ (uniform limits); so $C_b$ is closed and not dense (consider $\1_{[0,1/2]}$). The measure $w$ is submodular (if $A\cap B\ne\emptyset$ both sides of the submodularity inequality equal $2$; otherwise $w(A\cup B)\le w(A)+w(B)$), countably subadditive and continuous from below, and it has \OS, \Sone{} and \SK; it is neither order continuous nor regular, in accordance with \cref{cor:density-Lusin}.
\end{example}

\begin{example}[Distorted probabilities]\label{ex:distorted}
Let $P$ be a Borel probability measure on a metric space $X$, $g\colon[0,1]\to[0,\infty)$ continuous and strictly increasing with $g(0)=0$, and $\mu:=g\circ P$ on $\cB$. Then
\begin{itemize}[leftmargin=1em]
\item $\mu$ satisfies \UA{}. Indeed, $g$ is uniformly continuous on $[0,1]$ with modulus $\omega$, and $g^{-1}$ is continuous at $0$. If $\mu(B)<\delta$, then $P(B)\le g^{-1}(\delta)=:\eta$, so $\mu(A\cup B)=g(P(A\cup B))\le g(\min(1,P(A)+\eta))\le g(P(A))+\omega(\eta)$, and $\omega(g^{-1}(\delta))\to0$ as $\delta\to0$.
\item $\mu$ has \Sone: if $\mu(A_n)\to0$, then $P(A_n)\to0$; a subsequence satisfies $\sum_iP(A_{n_i})<\infty$, so $P(\bigcup_{i\ge k}A_{n_i})\to0$ and hence $\mu(\bigcup_{i\ge k}A_{n_i})\to0$.
\item $\mu$ is regular, since $P$ is regular on metric spaces (approximation by closed sets from inside and by open sets from outside) and $\mu(G\setminus F)=g(P(G\setminus F))$. Hence $\mu$ has \OS{} (\cref{lem:implications}(f)) and \Lu$[\mu]$ (the Lusin sets of $P$ work). By \cref{cor:density-Lusin} and \cref{thm:closure}, $C_b(X)\cap\Lone(\mu)$ is dense in $\Lone(g\circ P)$, and every $f\in\Lone(\mu)$ is a Choquet mean limit of functions in $C_b(X)$.
\item For $g(s)=s^{1/p}$, $1<p<\infty$, we have $\nC f=\int_0^\infty P(\{|f|>t\})^{1/p}\,dt=\frac1p\|f\|_{L^{p,1}(P)}$, the Lorentz norm with the normalization $\|f\|_{p,1}=p\int_0^\infty P(\{|f|>t\})^{1/p}\,dt$. This is consistent with the classical density of continuous functions in $L^{p,1}$.
\item Strict monotonicity of $g$ is used only to deduce that $P(B)$ is small when $\mu(B)$ is small; we make no claim when $g$ vanishes on an interval.
\end{itemize}
\end{example}

\begin{example}[A measure which is not null-additive]\label{ex:independence}
Let $X=[0,1]$, $P$ the Lebesgue measure and $\mu(A)=P(A)(1+\1_A(0))$ on $\cB$. By the proof of \cref{lem:implications}(d), $\mu$ is $2$-relaxed subadditive, hence satisfies \SK, but it is not null-additive. Thus \cref{thm:QC-approx}(b) and \cref{thm:Choquet-density} apply to $\mu$, while \cref{thm:QC-approx}(c) does not. The distorted probability $e^{-1/P}$, on the other hand, satisfies \UA{} but is not relaxed subadditive (\cref{lem:implications}(d)).
\end{example}
\noindent \Cref{tab:examples} summarizes the examples; the verification of the less obvious entries follows the table.

\begin{table}[H]
\centering
\small
\caption{Properties of the examples ($+$ holds, $-$ fails). ``Dense'' means that $C_b\cap\Lone$ is dense in $\Lone$.}\label{tab:examples}
\setlength{\tabcolsep}{4pt}
\begin{tabular}{@{}lcccccccccc@{}}
\toprule
Example & p.g.p. & \SK & \Sone & \PS & \OS & regular & order cont. & \Lu & dense \\
\midrule
\ref{ex:countable} (uncountable sets) & $+$ & $+$ & $+$ & $+$ & $-$ & $-$ & $-$ & $-$ & $-$ \\
\ref{ex:dense} (dense sets) & $-$ & $-$ & $-$ & $-$ & $+$ & $-$ & $-$ & $-$ & $-$ \\
\ref{ex:lambda} ($\lambda$) & $+$ & $+$ & $+$ & $+$ & $-$ & $-$ & $-$ & $-$ & $-$ \\
\ref{rem:S1-not-nec} ($v$) & $-$ & $-$ & $-$ & $-$ & $+$ & $+$ & $+$ & $+$ & $+$ \\
\ref{ex:nonempty} ($w$) & $+$ & $+$ & $+$ & $+$ & $+$ & $-$ & $-$ & $-$ & $-$ \\
\ref{ex:possibility} ($\Pi$) & $+$ & $+$ & $+$ & $+$ & $+$ & (*) & not claimed & (*) & (*) \\
\ref{ex:distorted} ($g\circ P$) & $+$ & $+$ & $+$ & $+$ & $+$ & $+$ & $+$ & $+$ & $+$ \\
$\Hd$, $\delta<n$, on $2^{\R^n}$ & $+$ & $+$ & $+$ & $+$ & $+$ & $-$ & $-$ & $-$ & $-$ \\
$\Hn$ on $\cL$ & $+$ & $+$ & $+$ & $+$ & $+$ & $+$ (**) & $-$ & & $+$ \\
\bottomrule
\end{tabular}

\smallskip
{\footnotesize (*) if and only if every $P_\varepsilon$ is discrete. (**) for sets of finite content (\cref{cor:delta-n}).}
\end{table}
\noindent \emph{Verification of some entries.} For \cref{ex:dense}: $\mu$ is not order continuous, since for a countable dense set $D=\{d_j\}$ the sets $D\setminus\{d_1,\dots,d_n\}$ are dense and decrease to $\emptyset$. It is not regular: for $A=\Q\cap[0,1]$, a closed $F\subseteq A$ is countable and hence nowhere dense, an open $G\supseteq A$ is dense, so $G\setminus F$ is dense. Continuous functions are not dense: $A=\Q\cap[0,1]$ (with $\mu(A)=1<\infty$) does not have property \RD. Indeed, if $\mu(E)<1$, then $E$ is not dense and misses a nondegenerate open interval $I\subseteq[0,1]$; from $A\setminus F\subseteq E$ we get $\Q\cap I\subseteq F$, so $I\subseteq F\subseteq G$ ($F$ is closed); then $I\setminus\Q\subseteq G\setminus A\subseteq E$, contradicting $I\cap E=\emptyset$. By \cref{thm:Choquet-density} (the implication (a)$\Rightarrow$(c) holds for every $\mu$), $C_b$ is not dense. For \cref{ex:lambda}: $[0,1/2]$ is not regular, since an open $G\supseteq[0,1/2]$ contains $[0,1/2+\eta)$, so $m(G\setminus F)\ge\eta$ for $F\subseteq[0,1/2]$. For \cref{rem:S1-not-nec}: $v$ is continuous from above (if $A_n\downto A\ne X$, pick $x\notin A$; then $x\notin A_n$ for large $n$), hence order continuous, and regular (for $A\ne X$ take $x\notin A$, $F=\emptyset$ and $G=X\setminus\{x\}$). For \cref{ex:distorted}: order continuity follows from the $\sigma$-additivity of $P$ and the continuity of $g$ at $0$. For $\Hd$: it is not order continuous, since $\Hd(\R^n\setminus B(0,k))=\infty$; and \Lu{} fails, since \Lu{} implies regularity (\cref{lem:indicator-QC}) while the Borel set $K$ of \cref{thm:cantor} is not regular.

%==========================================================================
\section{Concluding remarks and open problems}\label{sec:concl}

We have studied the approximation of measurable functions by continuous functions for monotone measures on metric spaces, in measure and with respect to the Choquet, Sugeno and Shilkret integrals. The classical description of the closure of the continuous functions in capacitary $L^1$-spaces as the class of quasicontinuous functions extends to monotone measures under property \Sone{} and outer smallness \OS{} (and \SK{} for unbounded functions in the Choquet case), and each of these conditions is needed. At the level of sets, the relevant notion is the symmetric approximation property \RD, which coincides with regularity under \OS. As applications we obtained corrected forms of mean approximation theorems of Cao and of Li and Yasuda, a description of the closure for possibility measures, and a Cantor-type example for the Hausdorff content. We hope that these observations may be useful to researchers working on non-additive measures and on capacitary function spaces.

Several questions remain open.
\begin{enumerate}[label=\textbf{Q\arabic*.},ref=Q\arabic*,leftmargin=3em]
\item\label{Q1} Let $\mu$ be a finite, continuous and weakly null-additive monotone measure which is \emph{not} null-additive, and let $f$ be \emph{unbounded} with $\Ch|f|\,d\mu<\infty$. Is $f$ a Choquet mean limit of bounded continuous functions? Is it, when $\mu(X\setminus U)=0$ with $U$ as in \cref{sec:Cao}, a Choquet mean limit of functions in $C_c(X)$, or, on $[a,b]$, of step functions? This is the only remaining gap in \cite[Theorem~4.2]{Cao15} on spaces with $\mu(X\setminus U)=0$ (in particular on locally compact spaces) and in \cite[Theorem~9]{LY04}; the bounded case and the null-additive case are settled by \cref{thm:Cao-char,thm:Cao-corrected}. The truncation method used above requires control of $\mu(\{|f|>M+t\}\cup N)$ for $\mu$-null sets $N$. We do not know a counterexample.
\item\label{Q2} Can property \Sone{} be replaced by property \PS{} in \cref{thm:approx-QC} under the p.g.p.\ and \OS?
\item\label{Q3} Is there a finite submeasure on a compact metric space such that $C_b$ is dense in $\Lone(\mu)$ but \Lu$[\mu]$ fails? By \cref{cor:density-Lusin}, such a measure must fail \Sone{} or \OS. Without \OS, property \RD{} is weaker than regularity: on $X=\{0\}\cup\{1/k:k\in\N\}$ let $\mu(A)=1$ if $A$ contains infinitely many points $1/k$, and $\mu(A)=0$ otherwise; then $\{0\}$ has property \RD{} (with $E=\{0\}$, $F=G=\emptyset$) but is not regular. In this example, however, $C_b$ is not dense, since $\{1/(2k):k\in\N\}$ does not have property \RD.
\item\label{Q4} Characterize the equality $\clM(C(X))=\M$ without order continuity. By \cref{thm:measure-density}(c), order continuity is not necessary, and the Lusin property is sufficient.
\item\label{Q5} Do \cref{thm:closure,thm:Choquet-density} extend to the Choquet--Lorentz spaces of \cite{KY23} and to Sugeno--Lorentz type spaces? For instance, for the functional $\Ch|f|^q\,d\mu^{q/p}$ one would need to know whether $\mu^{q/p}$ inherits \SK, \Sone{} and \OS{} from $\mu$; we have not checked this.
\end{enumerate}

%==========================================================================
\appendix
\section{Proofs of auxiliary results}\label{app:proofs}

\subsection{Proof of \texorpdfstring{\cref{lem:elementary}}{Lemma 2.5}}
\noindent (a) Put $\Phi(t):=\mu(\{g>t\})$ and $\Psi(t):=\mu(\{g\ge t\})$. Both are non-increasing, and $\Phi(t)\le\Psi(t)\le\Phi(s)$ for $s<t$. Let $t_0:=\sup\{t:\Phi(t)=\infty\}$ ($t_0:=0$ if this set is empty). On $[0,t_0)$ both functions are $\infty$. On $(t_0,\infty)$ the function $\Phi$ is finite and non-increasing, hence continuous except at countably many points; at a point $t$ of continuity, $\Psi(t)$ lies between $\Phi(t)$ and $\lim_{s\uparrow t}\Phi(s)=\Phi(t)$. So $\Phi=\Psi$ almost everywhere, and the Choquet integrals agree. For the Sugeno integral, ``$\le$'' is clear since $\{g>t\}\subseteq\{g\ge t\}$; conversely, for $t'<t$ we have $\min(t',\mu(\{g>t'\}))\ge\min(t',\mu(\{g\ge t\}))$, and letting $t'\uparrow t$ gives $\sup\ge\min(t,\mu(\{g\ge t\}))$. The Shilkret case is identical, with $t'\mu(\{g\ge t\})\to t\mu(\{g\ge t\})$ (if $\mu(\{g\ge t\})=\infty$, both sides are $\infty$).\\
\noindent (b) For $0<t<\alpha$ we have $\{g>t\}\supseteq\{g\ge\alpha\}$, so $\Ch g\,d\mu\ge\int_0^\alpha\mu(\{g>t\})\,dt\ge\alpha\mu(\{g\ge\alpha\})$. Also $t\mu(\{g\ge t\})\le\int_0^t\mu(\{g>s\})\,ds\le\Ch g\,d\mu$ for each $t$, so $\Sh g\,d\mu\le\Ch g\,d\mu$. The Sugeno bound is the definition.\\
\noindent (c) For $t\le\alpha$, $\min(t,\cdot)\le\alpha$; for $t>\alpha$, $\min(t,\mu(\{g\ge t\}))\le\mu(\{g\ge\alpha\})$.\\
\noindent (d) For $t>0$, $(g-M)^+>t$ if and only if $g>M+t$; substitute $s=M+t$.\\
\noindent (e) $\{c\1_A>t\}=A$ for $0\le t<c$ and $=\emptyset$ for $t\ge c$.\\
\noindent (f) $\{|cf|>t\}=\{|f|>t/|c|\}$ for $c\ne0$; substitute.\\
\noindent (g) $\Phi(t)=\mu(\{|f|>t\})$ is non-increasing with $\int_0^\infty\Phi<\infty$, hence finite for $t>0$ (if $\Phi(t_1)=\infty$, then $\Phi=\infty$ on $[0,t_1]$). $T(M)\to0$ by dominated convergence, and the last inequality is (b).\\
\noindent (h) By induction: $\mu(A_1\cup B)\le K(\mu(A_1)+\mu(B))$ with $B=A_2\cup\dots\cup A_m$ and $\mu(B)\le K^{m-2}\sum_{i\ge2}\mu(A_i)$. Then $\{|\sum_ig_i|>t\}\subseteq\bigcup_i\{|g_i|>t/m\}$, so $\nC{\sum_ig_i}\le K^{m-1}\sum_i\int_0^\infty\mu(\{|g_i|>t/m\})\,dt=mK^{m-1}\sum_i\nC{g_i}$. For the converse take $f=\1_A$ and $g=\1_{B\setminus A}$ and use (e).\\
\noindent (i) By (b), $\mu(\{|f_n-f|\ge\alpha\})\le\alpha^{-1}\Sh|f_n-f|\,d\mu\le\alpha^{-1}\Ch|f_n-f|\,d\mu$. Sugeno $\Rightarrow$ measure: once $\nS{f_n-f}<\alpha$, (b) gives $\min(\alpha,\mu(\{|f_n-f|\ge\alpha\}))<\alpha$, hence $\mu(\{|f_n-f|\ge\alpha\})\le\nS{f_n-f}\to0$. Measure $\Rightarrow$ Sugeno follows from (c). \qed

\subsection{Proof of \texorpdfstring{\cref{lem:implications}}{Lemma 2.6}}
\noindent (a) Given $\mu(A_n)\to0$, pick $n_i$ increasing with $\mu(A_{n_i})<2^{-i}$; then $\mu(\bigcup_{i\ge k}A_{n_i})\le2^{1-k}$. For \Sone$\Rightarrow$\PS, note that $\limsup_iA_{n_i}\subseteq\bigcup_{i\ge k}A_{n_i}$ for each $k$.\\
\noindent (b) Let $\delta(\cdot)$ be the function of the p.g.p.\ in $\varepsilon$--$\delta$ form (\cref{rem:pgp}); we may assume $\delta(\varepsilon)\le\varepsilon$. Put $\varepsilon_1:=1$ and $\varepsilon_{k+1}:=\min(\delta(\varepsilon_k),\varepsilon_k/2)$. Given $\mu(A_n)\to0$, choose $n_1<n_2<\cdots$ with $\mu(A_{n_k})<\varepsilon_{k+1}$. We claim that $\mu(A_{n_k}\cup\dots\cup A_{n_{k+l}})<\varepsilon_k$ for all $k,l\ge0$, and prove it by induction on $l$, for all $k$ simultaneously. For $l=0$: $\mu(A_{n_k})<\varepsilon_{k+1}\le\varepsilon_k$. For the induction step, the union is $A_{n_k}\cup U$ with $U=A_{n_{k+1}}\cup\dots\cup A_{n_{k+l}}$; by the induction hypothesis $\mu(U)<\varepsilon_{k+1}\le\delta(\varepsilon_k)$, and $\mu(A_{n_k})<\varepsilon_{k+1}\le\delta(\varepsilon_k)$; so by the p.g.p.\ the union has measure less than $\varepsilon_k$. Since $\varepsilon_k\le2^{1-k}\to0$, property \PC{} applied to $E_i:=A_{n_i}$ gives $\mu(\bigcup_{i\ge k}A_{n_i})\to0$. Under continuity from below one obtains directly $\mu(\bigcup_{i\ge k}A_{n_i})=\lim_l\mu(A_{n_k}\cup\dots\cup A_{n_{k+l}})\le\varepsilon_k$.\\
\noindent (c) For a submeasure, $\mu(A\cup B)\le\mu(A)+\mu(B)<\mu(A)+\varepsilon$ if $\mu(B)<\delta:=\varepsilon$. For a $K$-relaxed subadditive $\mu$, $\mu(A\cup B)\le K\mu(A)+K\mu(B)\le K\mu(A)+\varepsilon$ if $\mu(B)<\varepsilon/K$. \SK{} $\Rightarrow$ p.g.p.: if $\mu(A),\mu(B)<\min(\delta(\varepsilon/2),\varepsilon/(2K))$, then $\mu(A\cup B)\le K\mu(A)+\varepsilon/2<\varepsilon$. \UA{} $\Rightarrow$ null-additivity: if $\mu(B)=0$, then $\mu(A\cup B)\le\mu(A)+\varepsilon$ for all $\varepsilon>0$.\\
\noindent (d) Let $X=[0,1]$ and $P$ the Lebesgue measure. \emph{\UA{} but not relaxed subadditive:} let $g(s)=e^{-1/s}$ for $s>0$, $g(0)=0$, and $\mu:=g\circ P$. The function $g$ is continuous and strictly increasing on $[0,1]$, hence uniformly continuous with some modulus $\omega$, and $g^{-1}$ is continuous at $0$. If $\mu(B)<\delta$, then $P(B)\le g^{-1}(\delta)=:\eta$, so $\mu(A\cup B)=g(P(A\cup B))\le g(\min(1,P(A)+\eta))\le g(P(A))+\omega(\eta)$, and $\omega(g^{-1}(\delta))\to0$ as $\delta\to0$. But $\mu$ is not relaxed subadditive: for disjoint $A,B$ with $P(A)=P(B)=a$, $\mu(A\cup B)/(\mu(A)+\mu(B))=e^{-1/(2a)}/(2e^{-1/a})=e^{1/(2a)}/2\to\infty$ as $a\to0$. \emph{Relaxed subadditive ($K=2$) but not \UA:} let $\mu(A)=P(A)(1+\1_A(0))$. It is monotone: if $A\subseteq B$ and $0\notin A$, then $\mu(A)=P(A)\le P(B)\le\mu(B)$, and if $0\in A$, then $0\in B$ and $\mu(A)=2P(A)\le2P(B)=\mu(B)$. Moreover $\mu(A\cup B)\le2P(A\cup B)\le2(\mu(A)+\mu(B))$. But $\mu(\{0\})=0$ and, for $A=(0,1/2]$, $\mu(A\cup\{0\})=1\ne1/2=\mu(A)$; so $\mu$ is not null-additive, hence does not satisfy \UA{} by (c).\\
\noindent (e) Given $\mu(E)<\varepsilon$, pick an open $G\supseteq E$ with $\mu(G)<\varepsilon$.\\
\noindent (f) Let $\delta:=\delta_{\rm pgp}(\varepsilon)$. If $\mu(E)<\delta$, regularity gives $F\subseteq E\subseteq G$ with $\mu(G\setminus F)<\delta$. Since $G=F\cup(G\setminus F)$ and $\mu(F)\le\mu(E)<\delta$, we get $\mu(G)<\varepsilon$.\\
\noindent (g) For the first claim take $E=G\setminus F$: indeed $A\triangle F=A\setminus F\subseteq G\setminus F$ and $A\triangle G=G\setminus A\subseteq G\setminus F$. For the second, let $\varepsilon>0$ and $\delta:=\delta_{\rm OS}(\varepsilon)$; take $F\subseteq G$ and $E$ with $\mu(E)<\delta$ and $(A\triangle F)\cup(A\triangle G)\subseteq E$, and $O\in\cO$ with $O\supseteq E$ and $\mu(O)<\varepsilon$. Put $F':=F\setminus O\in\cC$ and $G':=G\cup O\in\cO$. (i) $F'\subseteq A$: if $x\in F\setminus O$, then $x\notin E$, so $x\notin A\triangle F$, and $x\in F$; hence $x\in A$. (ii) $A\subseteq G'$: if $x\in A\setminus O$, then $x\notin A\triangle G$, so $x\in G$. (iii) $G'\setminus F'\subseteq O$: if $x\in G'\setminus F'$ and $x\notin O$, then $x\in G$ and $x\notin F$; since $x\notin E$, $x\notin A\triangle G$ gives $x\in A$, and then $x\notin A\triangle F$ gives $x\in F$, a contradiction. So $\mu(G'\setminus F')\le\mu(O)<\varepsilon$. \qed

\subsection{Proof of \texorpdfstring{\cref{lem:compact-Lusin}}{Lemma 5.2}}\label{app:compact-Lusin}
\noindent Let $\rho$ be a compatible complete metric on $X$.\\
\noindent (1) \emph{Tightness.} Fix $\varepsilon>0$ and a dense sequence $(x_i)$. For $j\ge1$ put $U_{j,s}:=\bigcup_{i\le s}\bar B_\rho(x_i,1/j)$, so that $A^{(j)}_s:=X\setminus U_{j,s}\downto\emptyset$ as $s\to\infty$. We choose $s_1,s_2,\dots$ recursively. Let $B_0:=\emptyset$. Given $B_{j-1}$, the sets $B_{j-1}\cup A^{(j)}_s$ decrease (in $s$) to $B_{j-1}$. By continuity from above ($\mu$ is finite) we may choose $s_j$ with $\mu(B_{j-1}\cup A^{(j)}_{s_j})<\mu(B_{j-1})+\varepsilon2^{-j}$, and put $B_j:=B_{j-1}\cup A^{(j)}_{s_j}$. Then $\mu(B_j)<\varepsilon$ for all $j$, and $B_j\upto B_\infty:=\bigcup_jA^{(j)}_{s_j}$, so $\mu(B_\infty)\le\varepsilon$ by continuity from below. The set $K:=X\setminus B_\infty=\bigcap_jU_{j,s_j}$ is $\rho$-closed, hence complete, and it is totally bounded, since for each $j$ it is covered by $s_j$ balls of radius $2/j$. So $K$ is compact and $\mu(X\setminus K)\le\varepsilon$. Applying this with $\varepsilon=1/j$ we obtain compact sets $K'_j$; put $K_j:=K'_1\cup\dots\cup K'_j$. These compact sets increase, $D_j:=X\setminus K_j\downto D_\infty$, and $\mu(D_\infty)\le\mu(D_j)\le1/j$, so $\mu(D_\infty)=0$.\\
\noindent (2) \emph{Closed Lusin sets.} By \cite[Theorem~4]{LY04} (finite, continuous, weakly null-additive fuzzy measure on a metric space; real-valued Borel $f$) there are closed sets $F_k$ with $\mu(X\setminus F_k)<1/k$ such that $f|_{F_k}$ is continuous. Put $F'_k:=F_1\cup\dots\cup F_k$. Then $f|_{F'_k}$ is continuous by the pasting lemma, and $E_k:=X\setminus F'_k\downto E_\infty$ with $\mu(E_\infty)=0$.\\
\noindent (3) The sets $C_k:=F'_k\cap K_k$ are compact and increasing, $f|_{C_k}$ is continuous, and $X\setminus C_k=E_k\cup D_k\downto E_\infty\cup D_\infty=:N$ (both sequences decrease, so $\bigcap_k(E_k\cup D_k)=E_\infty\cup D_\infty$). By weak null-additivity $\mu(N)=0$, and by continuity from above $\mu(X\setminus C_k)\to\mu(N)=0$. \qed

%\section*{CRediT authorship contribution statement}
%[To be completed by the authors.]

%\section*{Declaration of competing interest}
%The authors declare that they have no known competing financial interests or personal relationships that could have appeared to influence the work reported in this paper.

%\section*{Data availability}
%No data was used for the research described in the article.

%\section*{Acknowledgements}

% AUTHORS: keep this statement only once you have actually done the verification it describes.
%\section*{AI-use disclosure}
%The authors utilized Claude (Opus 5.5) solely for minor grammatical refinements and formatting adjustments of their original draft. All conceptual development, mathematical proofs, and literature reviews were conducted entirely by the authors, who maintain full responsibility for the integrity and content of this paper.

%==========================================================================


\begin{thebibliography}{99}
\small

\bibitem{AH96}
D.R.~Adams, L.I.~Hedberg, \textit{Function Spaces and Potential Theory}, Grundlehren der mathematischen Wissenschaften, vol.~314, Springer, Berlin, 1996.

\bibitem{BCRS25}
R.~Basak, Y.-W.B.~Chen, P.~Roychowdhury, D.~Spector, The capacitary John--Nirenberg inequality revisited, arXiv:2501.11412, 2025.

\bibitem{Cao15}
S.~Cao, Regularity of fuzzy measures on complete and separable metric spaces, in: \textit{2015 12th International Conference on Fuzzy Systems and Knowledge Discovery (FSKD)}, IEEE, 2015, pp.~168--172.

\bibitem{CWL16}
S.~Cao, Y.~Wang, J.~Li, Approximation of fuzzy neural networks based on Choquet integral, \textit{J. Intell. Fuzzy Syst.} 31 (2016) 691--698.

\bibitem{CMS11}
J.~Cerd\`a, J.~Mart\'in, P.~Silvestre, Capacitary function spaces, \textit{Collect. Math.} 62 (2011) 95--118.

\bibitem{COS24}
Y.-W.B.~Chen, K.H.~Ooi, D.~Spector, Capacitary maximal inequalities and applications, \textit{J. Funct. Anal.} 286 (2024) 110396.

\bibitem{Den94}
D.~Denneberg, \textit{Non-Additive Measure and Integral}, Kluwer Academic Publishers, Dordrecht, 1994.

\bibitem{Gra14}
L.~Grafakos, \textit{Classical Fourier Analysis}, third ed., Graduate Texts in Mathematics, vol.~249, Springer, New York, 2014.

\bibitem{HH25}
P.~Harjulehto, R.~Hurri-Syrj\"anen, On Lebesgue points and measurability with Choquet integrals, \textit{J. Geom. Anal.} 36 (2026) 194. % Lưu ý: Kiểm tra lại nhãn HH25 và năm 2026 xem đã khớp chưa

\bibitem{JS96}
Q.~Jiang, H.~Suzuki, Fuzzy measures on metric spaces, \textit{Fuzzy Sets Syst.} 83 (1996) 99--106.

\bibitem{Kaw07}
J.~Kawabe, Regularity and Lusin's theorem for Riesz space-valued fuzzy measures, \textit{Fuzzy Sets Syst.} 158 (2007) 895--903.

\bibitem{Kec95}
A.S.~Kechris, \textit{Classical Descriptive Set Theory}, Springer, New York, 1995.

\bibitem{KY23}
J.~Kawabe, N.~Yamada, The completeness and separability of function spaces in nonadditive measure theory, \textit{Fuzzy Sets Syst.} 466 (2023) 108409.

\bibitem{LM11}
J.~Li, R.~Mesiar, Lusin's theorem on monotone measure spaces, \textit{Fuzzy Sets Syst.} 175 (2011) 75--86.

\bibitem{LMPK15}
J.~Li, R.~Mesiar, E.~Pap, E.P.~Klement, Convergence theorems for monotone measures, \textit{Fuzzy Sets Syst.} 281 (2015) 103--127.

\bibitem{LOM21}
J.~Li, Y.~Ouyang, R.~Mesiar, Generalized convergence theorems for monotone measures, \textit{Fuzzy Sets Syst.} 412 (2021) 53--64.

\bibitem{LY04}
J.~Li, M.~Yasuda, Lusin's theorem on fuzzy measure spaces, \textit{Fuzzy Sets Syst.} 146 (2004) 121--133.

\bibitem{Pap95}
E.~Pap, \textit{Null-Additive Set Functions}, Kluwer Academic Publishers, Dordrecht, 1995.

\bibitem{PS23}
A.C.~Ponce, D.~Spector, Some remarks on capacitary integrals and measure theory, in: S.~Lenhart, J.~Xiao (Eds.), \textit{Potentials and Partial Differential Equations: The Legacy of David R. Adams}, Advances in Analysis and Geometry, vol.~8, De Gruyter, Berlin, 2023.

\bibitem{ST22}
H.~Saito, H.~Tanaka, Dual of the Choquet spaces with general Hausdorff content, \textit{Studia Math.} 266 (2022) 323--335.

\bibitem{Sun94}
Q.~Sun, Property (S) of fuzzy measure and Riesz's theorem, \textit{Fuzzy Sets Syst.} 62 (1994) 117--119.

\bibitem{WK09}
Z.~Wang, G.J.~Klir, \textit{Generalized Measure Theory}, Springer, New York, 2009.

\bibitem{WYT13}
T.~Watanabe, S.~Yamada, T.~Tanaka, Note on the regularity of nonadditive measures, \textit{J. Appl. Math.} 2013 (2013) 234681.

\bibitem{WH94}
C.~Wu, M.~Ha, On the regularity of the fuzzy measure on metric fuzzy measure spaces, \textit{Fuzzy Sets Syst.} 66 (1994) 373--379.

\end{thebibliography}
\end{document}